# AN ERDŐS-KAC THEOREM FOR INTEGERS WITH DENSE DIVISORS

GÉRALD TENENBAUM AND ANDREAS WEINGARTNER

ABSTRACT. We show that for large integers $n$, whose ratios of consecutive divisors are bounded above by an arbitrary constant, the number of prime factors follows an approximate normal distribution, with mean $C \log_2 n$ and variance $V \log_2 n$, where $C = 1/(1 - \mathrm{e}^{-\gamma}) \approx 2.280$ and $V \approx 0.414$. This result is then generalized in two different directions.

## 1. Introduction

Let $\{d_j(n) : 1 \leqslant j \leqslant \tau(n)\}$ denote the increasing sequence of the divisors of a natural number $n$. Define

$$\mathcal{D}(t) := \{n \geqslant 1 : d_{j+1}(n) \leqslant t d_j(n) \quad (1 \leqslant j < \tau(n))\}, \quad \mathcal{D}(x,t) := \mathcal{D}(t) \cap [1,x],$$

and the cardinality $D(x,t) := |\mathcal{D}(x,t)|$. It was first proved in [10] that the arithmetic function $n \mapsto \max_j \log\{d_{j+1}(n)/d_j(n)\}/\log n$ $(n > 1)$ has a limiting distribution continuous on $[0,1]$ and strictly increasing on $[1/2,1]$. This implies that $D(x,t) = o(x)$ if $t = x^{o(1)}$ and $D(x,t) \leqslant (1-c_\varepsilon)x$ if $t \leqslant x^{1-\varepsilon}$. It is hence of clear arithmetical interest to describe the statistical structure of prime factors of integers satisfying such a strong multiplicative constraint. Note that [16, th. 1.3] provides sharp quantitative estimates for $D(x,t)$—see (5.2).

Theorem 1.1 establishes an Erdős-Kac theorem for integers in $\mathcal{D}(x,t)$. It shows that the number of prime factors of integers in $\mathcal{D}(x,t)$ follows a suitable Gaussian law, but with a mean that exceeds five times its variance. Therefore, the usual Poisson distribution for the number of prime factors does not apply here, and so the Gaussian distribution is not a consequence of a Poisson behavior.

Results of Erdős-Kac type usually rest on independence properties of prime factors of integers, as the historical work [3] and many subsequent developments, for instance, among others, [5]. In the situation studied here the constraints on divisors are so severe that they rule out any approach based on independence of prime factors: for instance, it is known (see, e.g., [2] and [14]) that large gaps must occasionally occur between prime factors of random integers, however this is strictly forbidden for the typical integers in the sets under study. This explains the difficulty and great technicality of the proofs, based on functional equations and Laplace inversion.

Let $\nu(n)$ denote the number of prime factors of an integer $n$, counted with or without multiplicity, i.e., with traditional notation, either $\nu(n) = \Omega(n)$ or $\nu(n) = \omega(n)$ for all $n$. Also let $y \mapsto \Phi(y)$ denote the standard normal distribution function and let $\log_k$ denote throughout the $k$-fold iterated logarithm. Finally, put

$$C := 1/(1 - \mathrm{e}^{-\gamma}) \approx 2.280291, \qquad V := C + 2K \approx 0.414284,$$

where $K$ is the constant defined in Lemma 4.4.

**Theorem 1.1.** *Let $t \geqslant 2$ be fixed. Uniformly for $x \geqslant 3$, $y \in \mathbb{R}$, we have*

$$\frac{|\{n \in \mathcal{D}(x,t) : \nu(n) \leqslant C \log_2 x + y\sqrt{V \log_2 x}\}|}{D(x,t)} = \Phi(y) + O_t\Big(\frac{1}{\sqrt{\log_2 x}}\Big).$$

We will extend this result in two different ways. The next statement describes the behavior of $\nu(n)$ on $\mathcal{D}(x,t)$, when $t$ is allowed to vary between 2 and $x$. For $x \geqslant 6$, $x \geqslant t \geqslant 2$, we define

$$\mu_{x,t} = C \log_2 x + (1-C)\log_2 t, \quad \sigma_{x,t}^2 := V \log_2 x + (1-V)\log_2 t.$$

**Theorem 1.2.** *Uniformly for $x \geqslant 6$, $x \geqslant t \geqslant 2$, $y \in \mathbb{R}$, we have*

$$\frac{|\{n \in \mathcal{D}(x,t) : \nu(n) \leqslant \mu_{x,t} + y\,\sigma_{x,t}\}|}{D(x,t)} = \Phi(y) + O\Big(\frac{1}{\sqrt{\log_2 x}}\Big). \tag{1.1}$$

Note that when $t = x$, Theorem 1.2 reduces to the classical Erdős-Kac theorem for natural numbers—see [13, th. III.4.15].

Theorem 1.2 adds to previous results in [19, th. 1], namely

$$\sum_{n \in \mathcal{D}(x,t)} \nu(n) = D(x,t)\{\mu_{x,t} + O(1)\},$$

$$\Big|\Big\{n \in \mathcal{D}(x,t) : |\nu(n) - \mu_{x,t}| > \xi\sqrt{\log_2 x}\Big\}\Big| \ll D(x,t)/\xi^2,$$

both estimates being valid uniformly for $x \geqslant t \geqslant 2$ and $\xi \geqslant 1$.

Next, we generalize Theorem 1.1 to a family of integer sequences that are defined as follows. Let $\vartheta$ be any arithmetic function and let $\mathcal{B} = \mathcal{B}_\vartheta$ be the set of positive integers containing $n = 1$ and all those $n \geqslant 2$ with prime factorization $n = p_1^{\alpha_1} \cdots p_k^{\alpha_k}$, $p_1 < p_2 < \ldots < p_k$, which satisfy

$$p_i \leqslant \vartheta\big(p_1^{\alpha_1} \cdots p_{i-1}^{\alpha_{i-1}}\big) \qquad (1 \leqslant i \leqslant k).$$

Let $\mathcal{B}(x) = \mathcal{B} \cap [1,x]$ and $B(x) = |\mathcal{B}(x)|$. When $\vartheta(n) = nt$, we have $\mathcal{B} = \mathcal{D}(t)$, by [11, lem. 2.2]. More generally, if $\vartheta(n)/n$ is a non-decreasing function of $n$, then

$$\mathcal{B} = \{n \geqslant 1 : d_{j+1}(n) \leqslant \vartheta(d_j(n)) \ \ (1 \leqslant j < \tau(n))\},$$

according to [17, th. 4].

When $\vartheta(n) = 1 + \sigma(n)$, where $\sigma(n)$ is the sum of the positive divisors of $n$, then $\mathcal{B}$ coincides with the set of *practical numbers* introduced by Srinivasan [8], i.e. integers $n$ such that every natural number $m \leqslant n$ can be expressed as a sum of distinct positive divisors of $n$. Stewart [9] and Sierpinski [7] showed that this is equivalent to the condition $p_i \leqslant 1 + \sigma(p_1^{\alpha_1} \cdots p_{i-1}^{\alpha_{i-1}})$ for $1 \leqslant i \leqslant k$.

**Theorem 1.3.** *Let $0 < a < 1$. Assume $\max(2,n) \leqslant \vartheta(n) \ll n\exp\{(\log n)^a\}$ for all $n \geqslant 1$. Uniformly for $x \geqslant 3$, $y \in \mathbb{R}$, we have*

$$\frac{|\{n \in \mathcal{B}(x) : \nu(n) \leqslant C \log_2 x + y\sqrt{V \log_2 x}\}|}{B(x)} = \Phi(y) + O\Big(\frac{1}{\sqrt{\log_2 x}}\Big). \tag{1.2}$$

Note that Theorem 1.3 implies Theorem 1.1. We will derive Theorem 1.3 from Theorem 10.1, which also leads to the following improvement of [19, cor. 1] regarding the average and normal order of $\nu(n)$ on $\mathcal{B}$.

**Theorem 1.4.** *Let $0 < a < 1$. Assume $\max(2,n) \leqslant \vartheta(n) \ll n\exp\{(\log n)^a\}$ for all $n \geqslant 1$. Uniformly for $x \geqslant 3$, $\xi \geqslant 1$, we have*

$$\sum_{n\in\mathcal{B}(x)} \nu(n) = B(x)\{C\log_2 x + O(1)\}, \tag{1.3}$$

$$\Big|\big\{n \in \mathcal{B}(x) : |\nu(n) - C\log_2 x| > \xi\sqrt{\log_2 x}\big\}\Big| \ll \frac{B(x)}{\xi^2}. \tag{1.4}$$

In Section 2, we consider sums of $z^{\nu(n)}$ $(z \in \mathbb{C})$ over integers $n$ without small prime factors, i.e. sifted integers. These estimates are then inserted into the functional equation in Lemma 3.1, which generalizes [16, Lemma 2.3]. This allows us to study sums of $z^{\nu(n)}$ over integers $n \in \mathcal{D}(x,t)$ (Section 6) and $n \in \mathcal{B}(x)$ (Section 10), with the help of Laplace transforms (Section 7). Finally, the Berry-Esseen inequality leads to the completion of the proof of Theorem 1.2 (Section 9) and Theorem 1.3 (Section 12). The proof of Theorem 1.4 is given in Section 11.

## 2. Summing $z^{\nu(n)}$ over sifted integers $n$

For $z \in \mathbb{C}$, let $u \mapsto \omega_z(u)$ be defined as the unique continuous function on $[1,\infty)$, satisfying the initial condition $u\omega_z(u) = z$ $(1 \leqslant u \leqslant 2)$ and the delay-differential equation

$$u\omega_z'(u) + \omega_z(u) = z\omega_z(u-1) \quad (u > 2). \tag{2.1}$$

Continuing $\omega_z$ by setting $\omega_z(u) = 0$ for $u < 1$, we see that (2.1) holds on $\mathbb{R} \smallsetminus \{1\}$. This function was introduced by Alladi [1], who established an Erdős-Kac theorem for sifted integers. It is a special case of a family of solutions of delay-differential equations studied in [4]. When $z = 1$, $\omega_1(u)$ is Buchstab's function, $\omega(u)$.

Define

$$I(s) := \int_0^s \frac{\mathrm{e}^t - 1}{t}\,\mathrm{d}t \qquad (s \in \mathbb{C}), \tag{2.2}$$

and let $\{b_k(z)\}_{k=0}^\infty$ denote the sequence of Taylor coefficients of $\mathrm{e}^{-zI(-s)}$ about $s = 0$. We thus have $b_0(z) = 1$, $b_1(z) = z$, $b_2(z) = \frac{1}{4}z(2z-1)$ and $b_k(z)$ is a polynomial in $z$ of degree $k$. Furthermore, put

$$a_k(z) := \sum_{0\leqslant j\leqslant k} \frac{(-1)^j}{j!} b_{k-j}(z) \qquad (k \geqslant 0),$$

so that $a_0(z) = 1$, $a_1(z) = z - 1$, $a_2(z) = \frac{1}{4}(z-2)(2z-1)$.

**Lemma 2.1.** *Let $K \geqslant 0$. For $z \in \mathbb{C}$, $1/2 \leqslant |z| \leqslant 2$ and $u \geqslant 1$, we have*

$$\omega_z(u) = \mathrm{e}^{-\gamma z} \sum_{0\leqslant k\leqslant K} \frac{b_k(z)u^{z-1-k}}{\Gamma(z-k)} + O_K\Big(u^{z-2-K}\Big), \tag{2.3}$$

$$\omega_z'(u) = \mathrm{e}^{-\gamma z} \sum_{0\leqslant k\leqslant K} \frac{b_k(z)u^{z-2-k}}{\Gamma(z-k-1)} + O_K\Big(u^{z-3-K}\Big). \tag{2.4}$$

*Moreover, for $u \geqslant 0$,*

$$\omega_z(u) = \mathrm{e}^{-\gamma z} \sum_{0\leqslant k\leqslant K} \frac{a_k(z)(u+1)^{z-1-k}}{\Gamma(z-k)} + O_K\Big((u+1)^{z-2-K}\Big), \tag{2.5}$$

$$\omega_z'(u) = \mathrm{e}^{-\gamma z} \sum_{0\leqslant k\leqslant K} \frac{a_k(z)(u+1)^{z-2-k}}{\Gamma(z-k-1)} + O_K\Big((u+1)^{z-3-K}\Big). \tag{2.6}$$

*We also have*

$$|\omega_z(u)| \leqslant |z|u^{|z|-1} \qquad (u>0,\ z\in\mathbb{C}). \tag{2.7}$$

*Proof.* The estimates (2.3) and (2.4) follow from [4, ths. 1 & 2]. The factor $\mathrm{e}^{-\gamma z}$ may be obtained by a computation similar to [12, (2.13)]—see also [1, th. 4]. The estimates (2.5) and (2.6) follow since for $\alpha\in\mathbb{C}$, $u\geqslant 1$,

$$u^\alpha = (u+1)^\alpha\left(1-\frac{1}{u+1}\right)^\alpha = (u+1)^\alpha\sum_{j\geqslant 0}\frac{(-1)^j}{(u+1)^j}\binom{\alpha}{j}.$$

Note that (2.5) and (2.6) are trivial when $0\leqslant u<1$, since $\omega_z(u)=\omega_z'(u)=0$ for $u<1$. The bound (2.7) coincides with [1, lemma 1]. □

Let $P^+(n)$—resp. $P^-(n)$—denote the largest—resp. the smallest—prime factor of an integer $n>1$ and make the convention that $P^+(1):=1$, $P^-(1):=\infty$. We define

$$\Phi_\nu(x,y,z) := \sum_{\substack{1\leqslant n\leqslant x\\ P^-(n)>y}} z^{\nu(n)} \qquad (x,y\in\mathbb{R},\ z\in\mathbb{C}). \tag{2.8}$$

**Theorem 2.2.** *Let $r\in[\frac12,2)$ be fixed. For $z\in\mathbb{C}$, $\frac12\leqslant|z|\leqslant r$, $2\leqslant y\leqslant x$, $u:=(\log x)/\log y$, we have*

$$\Phi_\nu(x,y,z) = \frac{x\omega_z(u)-zy}{\log y} - \frac{x\mathrm{e}^{-\gamma z}u^{z-2}}{\Gamma(z-1)(\log y)^2} + O_r\left(x\frac{u^{z-3}+E}{(\log y)^2}\right), \tag{2.9}$$

*where $E=(\log x)^{2|z|+1}\mathrm{e}^{-\sqrt{\log y}}$.*

*Proof.* This is an improved version of [1, th. 3]. Just like in [1], it follows from [1, th. 2] by inserting (2.4) with $K=0$ (instead of [1, (4.2)]) into [1, (4.3)]. Note that there is a typo in [1, (4.3)]: the + sign in front of the integral should be a − sign. □

Some notation is necessary for the next statement. We put

$$H_\nu(s,z) := \zeta(s)^{-z}\sum_{n\geqslant 1}\frac{z^{\nu(n)}}{n^s} = \begin{cases}\displaystyle\prod_p\left(1-\frac{z}{p^s}\right)^{-1}\left(1-\frac{1}{p^s}\right)^z & \text{if } \nu=\Omega,\\ \displaystyle\prod_p\left(1+\frac{z}{p^s-1}\right)\left(1-\frac{1}{p^s}\right)^z & \text{if } \nu=\omega,\end{cases} \tag{2.10}$$

where the infinite products converge uniformly when $\sigma=\operatorname{Re}s\geqslant\frac23$, $|z-1|\leqslant\frac12$, say. Furthermore, we write

$$G_\nu(s,y,z) := \left(\sum_{P^+(n)\leqslant y}\frac{z^{\nu(n)}}{n^s}\right)^{-1} = \begin{cases}\displaystyle\prod_{p\leqslant y}\left(1-\frac{z}{p^s}\right) & \text{if } \nu=\Omega,\\ \displaystyle\prod_{p\leqslant y}\left(1+\frac{z}{p^s-1}\right)^{-1} & \text{if } \nu=\omega.\end{cases}$$

Finally, we define

$$h_\nu(s,y,z) := \frac{H_\nu(s,z)G_\nu(s,yz)}{\Gamma(z)},\quad h_\nu(y,z):=h_\nu(1,y,z),$$

$$J_\nu(y,z) := \frac{H_\nu'(1,z)}{H_\nu(1,z)}+\frac{G_\nu'(1,z,y)}{G_\nu(1,z,y)}+\gamma z-1.$$

**Theorem 2.3.** *Let $r < 2$ be fixed. For $z \in \mathbb{C}$, $|z| \leqslant r$ and $\frac{3}{2} \leqslant y \leqslant x$, we have*

$$\begin{aligned}\Phi_\nu(x,y,z) \\ = x(\log x)^{z-1}\left\{h_\nu(y,z)\left(1+\frac{z-1}{\log x}J_\nu(y,z)\right)+O\left(\frac{(\log y)^7}{(\log x)^2}\right)\right\}.\end{aligned} \tag{2.11}$$

*Proof.* We apply [13, th. II.5.4], with parameters $N = 1$, $w = \max(1 - \operatorname{Re} z, 0)$, hence $w \leqslant 3$, and select $g_z = g_{z,y}$ as the multiplicative function defined, for $p \leqslant y$, by

$$\sum_{k\geqslant 0} g_{z,y}(p^k)\xi^k = (1-\xi)^z \quad (|\xi| < 1),$$

and, for $p > y$, by

$$\sum_{k\geqslant 0} g_{z,y}(p^k)\xi^k = \begin{cases} \dfrac{(1-\xi)^z}{1-z\xi} & \text{if } \nu = \Omega,\ |\xi| < \min(1, 1/|z|), \\ (1-\xi)^z\left(1+\dfrac{z\xi}{1-\xi}\right) & \text{if } \nu = \omega,\ |\xi| < 1. \end{cases}$$

By [13, th. II.5.4], it thus suffices to show that

$$E(z,y) := \sum_{n\geqslant 1} |g_{z,y}(n)|\frac{(\log 3n)^5}{n} \ll (\log y)^7.$$

When $p \leqslant y$, we have $g_{z,y}(p) = -z$ and Cauchy's formula yields $|g_{z,y}(p^k)| \leqslant 4$ for $k \geqslant 2$ since $|z| \leqslant r < 2$. When $p > y$, we have $g_{z,y}(p) = 0$ and, following [13, pp. 300-301], we have $|g_{z,y}(p^k)| \leqslant M(2-\delta)^k$ for $k \geqslant 2$, where $M$ and $\delta > 0$ depend on $r$ only.

Since $\log(3mn) \leqslant \log(3m)\log(3n)$ for all $m, n \geqslant 1$, we may write

$$E(z,y) \leqslant \sum_{P^+(n)\leqslant y} |g_{z,y}(n)|\frac{(\log 3n)^5}{n} \sum_{P^-(n)>y} |g_{z,y}(n)|\frac{(\log 3n)^5}{n} =: E_1E_2,$$

say. Now

$$E_2 \leqslant \prod_{p>y}\left(1+\sum_{k\geqslant 2}\frac{M(2-\delta)^k(\log 3p^k)^5}{p^k}\right) \ll 1.$$

To estimate $E_1$, we define

$$F_z(s) := \prod_{p\leqslant y}\left(1+\frac{|z|}{p^s}+\sum_{k\geqslant 2}\frac{4}{p^{ks}}\right).$$

By Cauchy's formula,

$$|F_z^{(5)}(1)| \leqslant \frac{5!}{\eta^5}|F_z(1-\eta)|,$$

for $0 < \eta < 1$. With $\eta = 1/(A\log y)$ we obtain

$$E_1 \leqslant |F_z^{(5)}(1)| \ll (A\log y)^5(\log y)^{|z|e^{1/A}} \ll (\log y)^7,$$

where $A := 1/\log(2/r)$. □

Selecting $y = 3/2$ in Theorem 2.3 yields the following statement.

**Lemma 2.4.** *Let $r < 2$. For $|z| \leqslant r$ we have*

$$\sum_{n\leqslant x} z^{\nu(n)} = x(\log x)^{z-1}\Big\{\lambda_{\nu,0}(z) + \frac{\lambda_{\nu,1}(z)}{\log x} + O_r\Big(\frac{1}{(\log x)^2}\Big)\Big\},$$

*where* $\lambda_{\nu,0}(z) := h_\nu(1,z),\ \lambda_{\nu,1}(z) := h_\nu(1,z)J_\nu(1,z)(z-1)$.

**Lemma 2.5.** *Let $r < 2$. Uniformly for $z \in \mathbb{C}$, $|z| \leqslant r$ and $y \geqslant 3/2$, we have*

$$J_\nu(y,z) = z\log y - 1 + O\big(\mathrm{e}^{-\sqrt{\log y}}\big), \tag{2.12}$$

$$h_\nu(y,z) = \frac{\mathrm{e}^{-\gamma z}}{\Gamma(z)(\log y)^z} + O\big(\mathrm{e}^{-\sqrt{\log y}}\big). \tag{2.13}$$

*Proof.* The first estimate follows by partial summation from the prime number theorem in the form

$$\sum_{p\leqslant y} \frac{\log p}{p-1} = \log y - \gamma + O\big(\mathrm{e}^{-\sqrt{\log y}}\big).$$

The second estimate is a consequence of a strong form of Mertens' product formula. □

## 3. The functional equation

We consider the following hypotheses regarding an arithmetic function $\vartheta$:

$$\vartheta : \mathbb{N}^* \to \mathbb{R} \cup \{\infty\}, \quad \vartheta(1) \geqslant 2, \quad \vartheta(n) \geqslant P^+(n) \quad (n \geqslant 2). \tag{3.1}$$

Here and throughout we let $\mathbb{N} := \{0, 1, 2, \cdots\}$ denote the set of non-negative integers, and define $\mathbb{N}^* := \mathbb{N} \smallsetminus \{0\}$. Recall notation (2.8).

**Lemma 3.1.** *Assume* (3.1). *For $z \in \mathbb{C}$ and $x \geqslant 0$ we have*

$$\sum_{m\leqslant x} z^{\nu(m)} = \sum_{n\in\mathcal{B}} z^{\nu(n)}\Phi_\nu\Big(\frac{x}{n}, \vartheta(n), z\Big).$$

*Proof.* This is an immediate consequence of the fact that each integer $m \geqslant 1$ factors uniquely as $m = nr$, with $n \in \mathcal{B}$, $P^-(r) > \vartheta(n)$. □

The following result is a consequence of [20, th. 4].

**Lemma 3.2.** *Let $0 \leqslant a < 1$ and assume* $\max(2,n) \leqslant \vartheta(n) \ll n\exp\{(\log n)^a\}$ *for all $n \geqslant 1$. There is a constant $c_\vartheta > 0$ such that*

$$B(x) = \frac{c_\vartheta x}{\log x}\Big\{1 + O\Big(\frac{1}{(\log x)^{1-a}}\Big)\Big\}.$$

**Lemma 3.3.** *Let $0 \leqslant a < 1$ and assume* $\max(2,n) \leqslant \vartheta(n) \ll n\exp\{(\log n)^a\}$ *for all $n \geqslant 1$. Let $z \in \mathbb{C}$ with $|z| = 1$,* $\operatorname{Re} z \geqslant 1/2$. *We have*

$$\lambda_{\nu,0}(z) = \sum_{n\in\mathcal{B}} \frac{z^{\nu(n)}}{n} h_\nu(\vartheta(n), z). \tag{3.2}$$

*Proof.* Lemma 3.1 yields

$$\sum_{m\leqslant x} z^{\nu(n)} = O(B(x)) + \sum_{n\in\mathcal{B}(x)} z^{\nu(n)}\left\{\Phi_\nu(x/n, \vartheta(n), z) - 1\right\}.$$

We divide both sides by $x(\log x)^{z-1}$ and estimate the left-hand side with Lemma 2.4 to get

$$\lambda_{\nu,0}(z)+o(1)=\frac{O(B(x))}{x(\log x)^{z-1}}+\sum_{n\geqslant 1} z^{\nu(n)}\mathbf{1}_{\mathcal{B}(x)}(n)\frac{\Phi_\nu(x/n,\vartheta(n),z)-1}{x(\log x)^{z-1}}, \tag{3.3}$$

where $\mathbf{1}_A(n)$ denotes the characteristic function of the set $A$. The modulus of the $n$-th term of the last series is bounded above, uniformly for $x\geqslant 1$, by

$$b_z(n):=\frac{\mathbf{1}_{\mathcal{B}}(n)}{n(\log\vartheta(n))^{\operatorname{Re} z}},$$

according to Theorems 2.2 and 2.3. We have $B(x)\asymp x/\log x$, by Lemma 3.2. Since $\operatorname{Re} z\geqslant 1/2$, the series $\sum_{n\geqslant 1} b_z(n)$ converges and

$$\lim_{x\to\infty}\frac{B(x)}{x(\log x)^{z-1}}=0.$$

Letting $x\to\infty$ in (3.3), the dominated convergence theorem allows interchange of limit and summation. Theorem 2.3 now yields the desired result. □

## 4. Zeros and poles

Define $M_z(u):=\omega_z(\mathrm{e}^u-1)$ $(u\geqslant 0)$ and

$$Q_z(s):=\int_0^\infty u^s\Big(\mathrm{e}^{zJ(u)}-1\Big)\,\mathrm{d}u \qquad (\operatorname{Re} s>\operatorname{Re} z-1),$$

where

$$J(u):=\int_u^\infty \frac{\mathrm{e}^{-t}}{t}\,\mathrm{d}t \qquad (u\geqslant 0).$$

From the definition of $\omega_z(u)$ we find that $\widehat{\omega}_z(s)=\mathrm{e}^{zJ(s)}-1$ by a computation similar to that of the proof of [13, th. III.6.7], corresponding to the case $z=1$. Now, we know from [13, lemma III.5.9] that, with notation (2.2), we have

$$u\mathrm{e}^{J(u)}=\mathrm{e}^{-\gamma-I(-u)} \quad (u>0). \tag{4.1}$$

Therefore

$$\mathrm{e}^{zJ(u)}=\sum_{k\geqslant 0} c_k(z)u^{k-z} \quad (u>0),$$

where $\{c_k(z)\}_{k=0}^\infty$ is the sequence of Taylor coefficients of $\mathrm{e}^{-\gamma z-zI(-u)}$ at the origin, so that $c_k(z)=\mathrm{e}^{-\gamma z}b_k(z)$ $(k\geqslant 0)$. This implies

$$Q_z(s)=\frac{-1}{s+1}+\sum_{k\geqslant 0}\frac{c_k(z)}{k+1+s-z}+\int_1^\infty u^s\Big\{\mathrm{e}^{zJ(u)}-1\Big\}\,\mathrm{d}u, \tag{4.2}$$

where the last integral converges since $J(u)\sim \mathrm{e}^{-u}/u$ for $u\geqslant 1$. This defines a meromorphic continuation of $Q_z(s)$ to the whole complex plane whose singularities are simple poles at $-1$ and $z-k-1$ $(k\geqslant 0)$—with an obvious abuse of notation if $z=1$.

**Lemma 4.1.** *For $s,z\in\mathbb{C}$ we have $(s+1)Q_z(s)=z\Gamma(s+1)\big\{1+\widehat{M}_z(s)\big\}$.*

*Proof.* For $\operatorname{Re} s > \operatorname{Re} z - 1$, we have

$$\begin{aligned}
Q_z(s) &= \int_0^\infty u^s \widehat{\omega}_z(u)\,\mathrm{d}u = \int_0^\infty u^s \int_0^\infty \omega_z(v)\mathrm{e}^{-uv}\,\mathrm{d}v\,\mathrm{d}u \\
&= \int_0^\infty \omega_z(v) \int_0^\infty u^s \mathrm{e}^{-uv}\,\mathrm{d}u\,\mathrm{d}v \\
&= \int_0^\infty \omega_z(v) \int_0^\infty \Big(\frac{w}{v}\Big)^s \mathrm{e}^{-w}\frac{\mathrm{d}w}{v}\,\mathrm{d}v = \Gamma(s+1)\int_0^\infty v\omega_z(v)\frac{\mathrm{d}v}{v^{s+2}} \\
&= \frac{z\Gamma(s+1)}{s+1}\left(1 + \int_0^\infty \omega_z(v)\frac{\mathrm{d}v}{(v+1)^{s+1}}\right) = \frac{z\Gamma(s+1)}{s+1}\Big\{1 + \widehat{M}_z(s)\Big\},
\end{aligned}$$

where the penultimate equation follows from (2.1) and partial integration. Since $Q_z(s)$ extends to a meromorphic function on all of $\mathbb{C}$, so does $\widehat{M}_z(s)$. □

**Corollary 4.2.** *Assume $z \notin \mathbb{Z}$. Then $\widehat{M}_z(s)$ has a simple pole at every $s = z-k-1$ ($k \in \mathbb{N}$) for which $c_k(z) \neq 0$, and no other singularities. Furhermore, $1+\widehat{M}_z(s)$ has simple zeros at negative integers. All other zeros of $1+\widehat{M}_z(s)$ are zeros of $Q_z(s)$.*

**Corollary 4.3.** *For $z \in \mathbb{C}$ define*

$$f_z(s) := \int_0^\infty u^s \mathrm{e}^{-u+zJ(u)}\,\mathrm{d}u \qquad (\operatorname{Re} s > \operatorname{Re} z - 1). \tag{4.3}$$

*Then $f_z(s)$ extends meromorphically to the entire complex plane via the formula*

$$f_z(s) = \Gamma(s+1)\big\{1+\widehat{M}_z(s)\big\} = \frac{(s+1)Q_z(s)}{z}\cdot \tag{4.4}$$

*Proof.* Integration by parts yields $Q_z(s) = zf_z(s)/(s+1)$, which implies the result by (4.2) and Lemma 4.1. □

For the statement of our next lemma we use the notation

$$\begin{aligned}
&T(s) := \int_0^s \frac{1-\mathrm{e}^{-t}}{t}\,\mathrm{d}t, \quad \mathrm{e}^{-\gamma z + zT(s)} = \sum_{k\geqslant 0} c_k(z)s^k \qquad (s\in\mathbb{C}), \\
&A := \mathrm{e}^{-\gamma}\Big\{\int_0^1 \frac{\mathrm{e}^{T(u)}-1-u}{u^2}\,\mathrm{d}u - 1\Big\} \\
&W(s,z) := \int_1^\infty u^s(\mathrm{e}^{zJ(u)}-1)\,\mathrm{d}u, \quad W := W(-1,1), \\
&B := A+W, \quad C := \frac{1}{1-\mathrm{e}^{-\gamma}}, \quad K := \mathrm{e}^{-\gamma}C^2(1-\gamma+BC).
\end{aligned} \tag{4.5}$$

We note the numerical values $C \approx 2.280291$, $K \approx -0.933003$.

**Lemma 4.4.** *Let $f_z(s)$ be given by (4.4). For $|z-1|$ sufficiently small and suitable $r>0$, the equation $f_z(s)=0$ has a unique solution $s=s_0(z)$ in the disk $|s+1| \leqslant r$. Moreover*

$$s_0(z) = -1 + C(z-1) + K(z-1)^2 + O((z-1)^3). \tag{4.6}$$

*Proof.* By Rouché's theorem we may assert that the equation $f_z(s) = 0$ has a unique solution in a suitable disk $|s+1| \leqslant r$ provided $|z-1|$ is sufficiently small, the details being left to the reader.

Write $s=-1+h$, $z=1+w$, so that the equation can be rewritten as

$$\sum_{k\geqslant 0}\frac{c_k(z)h}{k-1+h-w}+hW(s,z)=1. \tag{4.7}$$

Since $c_1(z)=z\mathrm{e}^{-\gamma z}$, we derive from the above that

$$\frac{z\mathrm{e}^{-\gamma z}h}{h-w}=1+O(h),$$

and in turn $h(z\mathrm{e}^{-\gamma z}-1)=-w+O(h^2+|hw|)$. Hence $h\ll w$ and

$$h=\frac{w}{1-\mathrm{e}^{-\gamma}}+O\big(w^2\big)=Cw+O(w^2).$$

To determine the second order term, we carry back into (4.7). Since

$$\sum_{k\neq 1}\frac{c_k(1)}{k-1}=A,$$

we have

$$\begin{aligned}\frac{zh\mathrm{e}^{-\gamma z}}{h-w}&=\frac{h\mathrm{e}^{-\gamma}+hw\mathrm{e}^{-\gamma}(1-\gamma)+O(h^3)}{h-w},\\ \sum_{k\neq 1}\frac{c_k(z)h}{k-1+h-w}&=Ah+O(h^2),\\ hW(s,z)&=hW+O(h^2),\end{aligned}$$

from which it follows that

$$\begin{aligned}h\mathrm{e}^{-\gamma}\Big\{1+w(1-\gamma)\Big\}&=(h-w)\Big\{1-hB\Big\}+O(w^3)\\ &=h-w-\mathrm{e}^{-\gamma}BCwh+O(w^3),\\ h=\frac{-w+O(w^3)}{\mathrm{e}^{-\gamma}-1+w\mathrm{e}^{-\gamma}(1-\gamma+BC)}&=\frac{w+O(w^3)}{1-\mathrm{e}^{-\gamma}}\Big(1+w\frac{\mathrm{e}^{-\gamma}(1-\gamma+BC)}{1-\mathrm{e}^{-\gamma}}\Big)\\ &=Cw+Kw^2+O(w^3).\end{aligned}$$

□

**Lemma 4.5.** *For suitable $\delta>0$ and $|z-1|\leqslant\delta$ the equation $1+\widehat{M}_z(s)=0$ has exactly four roots, all simple, in the half-plane $\operatorname{Re}s\geqslant -151/50$. These roots are: $s_0(z)$, $-1$, $-2$, $-3$, where $s_0(z)$ is defined in Lemma 4.4.*

*Proof.* We write $\sigma=\operatorname{Re}s$, $\tau=\operatorname{Im}s$. For $|z|\leqslant 2$, $\sigma\geqslant 2$, the upper bound (2.7) implies

$$|\widehat{M}_z(s)|\leqslant\int_1^\infty\frac{|\omega_z(u)|}{(u+1)^{\sigma+1}}\,\mathrm{d}u\leqslant\int_1^\infty\frac{2u}{(u+1)^3}\,\mathrm{d}u=\frac{3}{4}\cdot$$

Thus, for $|z|\leqslant 2$, $1+\widehat{M}_z(s)$ does not vanish in the half-plane $\sigma\geqslant 2$.

The estimate (2.6) implies

$$\omega_z'(u)=\mathrm{e}^{-\gamma z}\sum_{0\leqslant k\leqslant 3}\frac{a_k(z)(u+1)^{z-2-k}}{\Gamma(z-k-1)}+O(u^{z-6})=:g_z(u)+O(u^{z-6}),$$

say. Integration by parts yields

$$\widehat{M}_z(s) = \frac{z}{s2^s} + \frac{1}{s}\int_1^\infty \frac{\omega_z'(u)}{(u+1)^s}\,\mathrm{d}u = \frac{z}{s2^s} + \frac{1}{s}\int_1^\infty \frac{\omega_z'(u) - g_z(u)}{(u+1)^s}\,\mathrm{d}u + \frac{h_z(s)}{s}, \tag{4.8}$$

where

$$h_z(s) := \int_1^\infty \frac{g_z(u)}{(u+1)^s}\,\mathrm{d}u = \sum_{0\leqslant k\leqslant 3} \frac{\mathrm{e}^{-\gamma z}2^{z-1-k-s}a_k(z)}{(s+k+1-z)\Gamma(z-k-1)}. \tag{4.9}$$

Since $(\omega_z'(u) - g_z(u)) \ll u^{z-6}$, the second integral in (4.8) converges for $\sigma > \operatorname{Re} z - 5$ and is bounded uniformly for $\sigma \geqslant \operatorname{Re} z - 249/50$, say. Thus, if $\operatorname{Re} z \leqslant 6/5$, (4.8) and (4.9) give the meromorphic extension of $\widehat{M}_z(s)$ to $\sigma \geqslant -7/2$ and show that $\widehat{M}_z(s) \ll 1/|\tau|$ for $\sigma \geqslant -7/2$ and $|\tau| \geqslant 1$, say. It follows that there exists a constant $\tau_0 > 0$ such that $1 + \widehat{M}_z(s)$ does not vanish if $|\tau| \geqslant \tau_0$ and $\sigma \geqslant -7/2$.

Integration by parts shows that the second integral in (4.8) is $\ll 1/s$, so that

$$\frac{1}{1+\widehat{M}_z(s)} = 1 - \frac{z}{s2^s} + O\Big(\frac{1}{\tau^2}\Big) \qquad (|\tau| \geqslant 1,\ \sigma \geqslant -7/2). \tag{4.10}$$

Because $1 + \widehat{M}_1(s)$ has a pair of complex zeros with $\operatorname{Re} s \approx -3.03 < -151/50$ (see [15, Prop. 1]), we define the rectangle $\mathcal{R}$ by $-151/50 \leqslant \sigma \leqslant 2$, $|\tau| \leqslant \tau_0$. As $z \to 1$, $\widehat{M}_z(s) \to \widehat{M}_1(s)$, uniformly on the boundary of $\mathcal{R}$. The argument principle guarantees that $1 + \widehat{M}_z(s)$ and $1 + \widehat{M}_1(s)$ have the same number of zeros, minus number of poles, inside $\mathcal{R}$, if $z$ is close enough to 1. However, according to [15, proofs of lem. 5 & prop. 1], inside $\mathcal{R}$, $1+\widehat{M}_1(s)$ has exactly one simple zero, at $s = -1$, and one simple pole, at $s = 0$. By Corollary 4.2, the only poles of $1+\widehat{M}_z(s)$ in $\mathcal{R}$ are at $s = z-1,\ z-2,\ z-3,\ z-4$, while $1+\widehat{M}_z(s)$ has zeros at $s = -1, -2, -3$, all of which are simple. It follows that $1 + \widehat{M}_z(s)$ must have exactly one other zero inside $\mathcal{R}$, which must be the zero defined in Lemma 4.4, since $f_z(s) = \Gamma(s+1)(1+\widehat{M}_z(s))$. □

## 5. The function $d_z(v)$

For $z \in \mathbb{C}$, define $d_z(v)$ for $v \in \mathbb{R}$ by $d_z(v) := 0$ if $v \leqslant 0$ and

$$d_z(v) := v^{z-1} - \int_0^{(v-1)/2} \frac{d_z(u)}{u+1}\,\omega_z\left(\frac{v-u}{u+1}\right)\mathrm{d}u \qquad (v > 0). \tag{5.1}$$

From [16, th. 1.3] we have

$$D(x,t) = x\,\eta_t\,d_1(v)\,\big\{1 + O\big(1/\log x\big)\big\} \qquad (x \geqslant t \geqslant 2,\ v := (\log x)/\log t), \tag{5.2}$$

where $\eta_t = 1 + O(1/\log t)$ and $d_1(v)$ satisfies [15, th. 1]

$$d_1(v) = \frac{C}{v+1}\Big\{1 + O\Big(\frac{1}{v^2}\Big)\Big\} \qquad (v \geqslant 1), \tag{5.3}$$

where $C$ is defined in (4.5). The next statement generalizes this to $d_z(v)$. We define

$$C_z := \frac{\Gamma(s_z+1-z)\Gamma(z)z}{((s_z+1)Q_z(s_z))'} \tag{5.4}$$

where $s_z := s_0(z)$ as defined in Lemma 4.4, and we note right away that it follows from Lemmas 4.1, 4.4, and (4.2), (4.6) that

$$C_z = C + O(z-1) \qquad (z \to 1). \tag{5.5}$$

**Theorem 5.1.** *For suitable* $\delta > 0$, $|z-1| \leqslant \delta$ *and* $v \geqslant 1$, *we have*

$$d_z(v) = C_z(v+1)^{s_0(z)}\big\{1 + O\big(1/v^2\big)\big\}.$$

*Proof.* We may assume $z \neq 1$, since the case $z = 1$ is covered by (5.3). Changing variables via

$$v = \mathrm{e}^w - 1, \quad q_z(w) = d_z(\mathrm{e}^w - 1), \quad M_z(w) = \omega_z(\mathrm{e}^w - 1),$$

equation (5.1) becomes

$$q_z(w) = (\mathrm{e}^w - 1)^{z-1} - \int_0^w q_z(u) M_z(w-u)\,\mathrm{d}u. \tag{5.6}$$

Since $|M_z(w)| \leqslant 2(e^w - 1)$ for $|z| \leqslant 2$ by (2.7), and $M_z(w) = 0$ for $w \leqslant \log 2$, it follows from (5.6) by induction on $\lfloor w/\log 2 \rfloor$ that $|q_z(w)| \ll \mathrm{e}^{2w}$ for $w \geqslant \log 2$, while $q_z(w) = (\mathrm{e}^w - 1)^{z-1}$ if $0 < w \leqslant \log 2$. Multiplying (5.6) by $\mathrm{e}^{-ws}$, where $\mathrm{Re}\, s > 2$, and integrating with respect to $w$ over $(0, \infty)$, we obtain the equation of Laplace transforms,

$$\widehat{q_z}(s) = \frac{\Gamma(s+1-z)\,\Gamma(z)}{\Gamma(s+1)} - \widehat{q_z}(s)\widehat{M_z}(s).$$

Indeed, the change of variables $u = \mathrm{e}^{-w}$ yields

$$\begin{aligned} \int_0^\infty (\mathrm{e}^w - 1)^{z-1}\mathrm{e}^{-ws}\,\mathrm{d}w &= \int_0^1 (1-u)^{z-1} u^{s-z}\,\mathrm{d}u \\ &= B(s-z+1, z) = \frac{\Gamma(s+1-z)(z)}{\Gamma(s+1)}. \end{aligned} \tag{5.7}$$

Solving for $\widehat{q_z}(s)$, we obtain

$$\widehat{q_z}(s) = \frac{\Gamma(s+1-z)\,\Gamma(z)}{\Gamma(s+1)\{1 + \widehat{M_z}(s)\}}.$$

By (4.4) and Lemma 4.4, $\widehat{q_z}(s)$ has a pole at $s = s_0(z)$. With the help of Lemmas 4.1 and 4.4, and equation (4.2), we find that the corresponding residue $C_z$ satisfies (5.4). From Corollary 4.2 we infer that the simple poles of $\Gamma(s+1-z)$ at $s = z-k$ ($k \geqslant 1$) are canceled by the simple poles of $\widehat{M_z}(s)$, while the zeros of $1 + \widehat{M_z}(s)$ at $s = -k$ ($k \geqslant 1$) are canceled by the poles of $\Gamma(s+1)$. It follows from Lemma 4.5 that the only pole of $\widehat{q_z}(s)$ with $\mathrm{Re}\, s \geqslant -151/50$ is at $s = s_0(z)$.

Inverting the Laplace transform, the residue theorem yields

$$q_z(w) = \frac{1}{2\pi i}\int_{3-i\infty}^{3+i\infty} \widehat{q_z}(s) e^{ws}\,\mathrm{d}s = C_z e^{s_0(z)w} + \frac{1}{2\pi i}\int_{\mathcal{C}} \widehat{q_z}(s) e^{ws}\,\mathrm{d}s,$$

where $\mathcal{C}$ is the contour of five line segments with endpoints

$$3 - i\infty,\ 3 - iT,\ -151/50 - iT,\ -151/50 + iT,\ 3 + iT,\ 3 + i\infty,$$

and $T = e^{8w}$.

We have

$$\log\Gamma(1+s-z) = \log\Gamma(s) - (z-1)\log s + O\big((z-1)/s\big),$$

since $(\log\Gamma(s))' = \log s + O(1/s)$ and $(\log\Gamma(s))'' \ll 1/s$. Together with (4.10), this yields

$$\widehat{q_z}(s) = \frac{\Gamma(z)}{s^z}\Big\{1 + O\Big(\frac{1}{s}\Big)\Big\},$$

on the contour $\mathcal{C}$. Thus, if $|z-1|\leqslant\delta$ and $\delta$ is sufficiently small,

$$\int_{\mathcal{C}}\widehat{q}_z(s)e^{ws}\,\mathrm{d}s\ll e^{-151w/50}\leqslant(v+1)^{s_0(z)-2},$$

by Lemma 4.4. □

## 6. Summing $z^{\nu(n)}$ over $\mathcal{D}(x,t)$

Put

$$D(x,t,z):=\sum_{n\in\mathcal{D}(x,t)}z^{\nu(n)},\qquad\Delta(x,t,z):=\frac{D(x,t,z)}{x}.\tag{6.1}$$

The following statement will be proved in Sections 7 and 8.

**Theorem 6.1.** *Let $\varepsilon>0$. For suitable $\delta>0$, the estimate*

$$\Delta(x,t,z)=\eta_{z,t}(\log t)^{z-1}d_z(v)\Big\{1+O\Big(\frac{1}{(\log x)^{1-\varepsilon}}\Big)\Big\},\tag{6.2}$$

*holds uniformly for $|z|=1$, $|z-1|\leqslant\delta$, $x\geqslant t\geqslant 2$, $v=(\log x)/\log t$. We have*

$$\eta_{z,t}=\eta_t\{1+O(z-1)\}=1+O\big(|z-1|+1/\log t\big),\tag{6.3}$$

*where $\eta_t$ is as in* (5.2).

We first derive some consequences.

**Corollary 6.2.** *For suitable $\delta>0$, $|z|=1$, $|z-1|\leqslant\delta$ and $v\geqslant 1$, we have*

$$\frac{d_z(v)}{d_1(v)}=v^{1+s_0(z)}\{1+O(z-1)\}.$$

Combining Theorem 6.1 with (5.2) and Corollary 6.2, we obtain the following statement.

**Corollary 6.3.** *Let $\varepsilon>0$. For suitable $\delta>0$, the estimate*

$$D(x,t,z)/D(x,t)=(\log x)^{1+s_0(z)}(\log t)^{z-2-s_0(z)}\left\{1+O(|z-1|+(\log x)^{\varepsilon-1})\right\},$$

*holds uniformly for $|z|=1$, $|z-1|\leqslant\delta$, $x\geqslant t\geqslant 2$.*

The proof of Corollary 6.2 assuming Theorem 6.1 rests upon the following estimates from [19, th. 1 & proof].

**Lemma 6.4.** *Let $\mu_{x,t}:=C\log_2 x+(1-C)\log_2 t$. For $x\geqslant 3$, $x\geqslant t\geqslant 2$, we have*

$$\sum_{n\in\mathcal{D}(x,t)}\nu(n)=D(x,t)\{\mu_{x,t}+O(1)\},\qquad\sum_{n\in\mathcal{D}(x,t)}\nu(n)^2\ll D(x,t)\mu_{x,t}^2.$$

*Proof of Corollary 6.2.* If $1/v^2\leqslant|z-1|$, the result follows from (5.3), (5.5), Theorem 5.1 and Lemma 4.4. Thus, we may assume $1/v^2\geqslant|z-1|$.

Writing $z^{\nu(n)}=\mathrm{e}^{i\varphi\nu(n)}=1+i\varphi\nu(n)+O(\varphi^2\nu(n)^2)$, Lemma 6.4 yields

$$\begin{aligned}D(x,t,z)/D(x,t)&=1+i\varphi\mu_{x,t}+O\big(|\varphi|+\varphi^2\mu_{x,t}^2\big)\\&=\mathrm{e}^{(z-1)\mu_{x,t}}\big\{1+O\big(|\varphi|+\varphi^2\mu_{x,t}^2\big)\big\},\end{aligned}\tag{6.4}$$

provided $\mu_{x,t}\varphi^2\ll 1$. By (6.2), (5.2) and (6.3), we obtain

$$\mathrm{e}^{(z-1)\mu_{x,t}}\big\{1+O\big(|\varphi|+\varphi^2\mu_{x,t}^2\big)\big\}=\frac{d_z(v)(\log t)^{z-1}}{d_1(v)}\Big\{1+O\Big(|z-1|+\frac{1}{(\log x)^{1-\varepsilon}}\Big)\Big\}.$$

We now choose $x$ such that $(\log_2 x)^2 = 1/|z-1|$. Note that $v^2 \leqslant 1/|z-1| = (\log_2 x)^2$ implies $t \geqslant 2$. We obtain

$$\begin{aligned}\frac{d_z(v)}{d_1(v)} &= \mathrm{e}^{(z-1)\mu_{x,t}}(\log t)^{1-z}\left\{1+O\left(z-1\right)\right\}\\ &= \Big(\frac{\log x}{\log t}\Big)^{C(z-1)}\left\{1+O\left(z-1\right)\right\}\\ &= v^{C(z-1)}\left\{1+O\left(z-1\right)\right\}.\end{aligned}$$

However, we have $C(z-1) = 1 + s_0(z) + O((z-1)^2)$ by Lemma 4.4 and

$$(\log v)(z-1)^2 \ll \mu_{x,t}\varphi^2.$$

This yields the required estimate. □

## 7. Proof of Theorem 6.1: auxiliary estimates

We shall use the notation

$$v_n := \frac{\log xt}{\log nt} \qquad (x \geqslant t \geqslant 2,\ n \geqslant 1), \qquad \alpha(z) := \frac{\mathrm{e}^{-\gamma z}}{\Gamma(z)} \qquad (z \in \mathbb{C}), \tag{7.1}$$

and, for $x \geqslant t \geqslant 2$, $|z| = 1$, $\mathrm{Re}\, z > 1/2$,

$$\delta(x,t,z) := \sum_{n\in\mathcal{D}(t)} \frac{z^{\nu(n)}\{\alpha(z)v_n^{z-1} - \omega_z\big(v_n-1\big)\}}{n\log nt}, \tag{7.2}$$

the series being convergent in view of (2.3). Here and throughout, we write $\log xt$ to mean $\log(xt)$. We also recall (6.1).

**Lemma 7.1.** *Let $\varepsilon > 0$. For $z \in \mathbb{C}$, $|z| = 1$, $\mathrm{Re}\, z > 1/2$ and $x \geqslant t \geqslant 2$, we have*

$$\begin{aligned}\Delta(x,t,z) = \delta(x,t,z) &+ \lambda_{\nu,0}(z)\{(\log x)^{z-1} - (\log xt)^{z-1}\} + c_{z,t}(\log xt)^{z-2}\\ &+ O\Big(\frac{\log t}{(\log x)^{2-\varepsilon}}\Big),\end{aligned} \tag{7.3}$$

*where $\lambda_{\nu,0}(z)$ is defined in Lemma 2.4 and $c_{z,t} \ll z-1$ only depends on $z$ and $t$.*

*Proof.* Put

$$X := \mathrm{e}^{(\log_2 x)^3}, \quad \mathcal{U}_1 := \{n \geqslant 1 : n^2 t \leqslant X\},$$
$$\mathcal{U}_2 := \{n \geqslant 1 : X < n^2 t \leqslant x\}, \quad \mathcal{U}_3 := \{n \geqslant 1 : \sqrt{x/t} < n \leqslant x\}$$

and define

$$S_j := \frac{1}{x}\sum_{\substack{n\in\mathcal{D}(x,t)\\ n\in\mathcal{U}_j}} z^{\nu(n)}\Phi_\nu\Big(\frac{x}{n}, nt, z\Big) \qquad (j = 1,2,3).$$

By Lemma 3.1 with $\vartheta(n) := nt$, we get

$$\frac{1}{x}\sum_{m\leqslant x} z^{\nu(m)} = S_1 + S_2 + S_3. \tag{7.4}$$

When $n \in \mathcal{U}_3$, we have $\Phi_\nu(x/n, nt, z) = 1$, so

$$S_3 = \Delta(x,t,z) + O\big(1/\sqrt{x}\big).$$

We then estimate $S_2$ with Theorem 2.2, and $S_1$ by Theorem 2.3.

Let us first show that the error terms from both theorems contribute an amount $\ll (\log t)/(\log xt)^{2-\varepsilon}$. Note that $n^2t \leqslant x$ implies $\log x/n \asymp \log x$.

In (2.9), the first error term contributes to $S_2$

$$\ll \sum_{\substack{n\in\mathcal{D}(t)\\ X<n^2t\leqslant x}} \frac{(\log x)^{z-3}}{n(\log nt)^{z-1}} \ll \sum_{n\in\mathcal{D}(x,t)} \frac{1}{n(\log x)^2} \ll \frac{(\log t)\log_2 x}{(\log x)^2} \ll \frac{\log t}{(\log x)^{2-\varepsilon}},$$

while the second error term yields

$$\ll \sum_{\substack{n\in\mathcal{D}(t)\\ X<n^2t\leqslant x}} \frac{(\log x)^3 \mathrm{e}^{-\sqrt{\frac{1}{2}\log X}}}{n} \ll \frac{1}{(\log x)^2},$$

since $nt > \sqrt{X}$ here.

The term $-zy/\log y$ in (2.9) contributes

$$\ll \frac{1}{x}\sum_{\substack{n\in\mathcal{D}(t)\\ nt\leqslant\sqrt{xt}}} \frac{nt}{\log nt} \ll \frac{\sqrt{xt}}{x\log xt}\sum_{\substack{n\in\mathcal{D}(t)\\ n\leqslant\sqrt{x/t}}} 1 \ll \frac{\log t}{(\log x)^2}.$$

The error term in (2.11) contributes to $S_3$

$$\ll \sum_{\substack{n\in\mathcal{D}(t)\\ n^2t\leqslant X}} \frac{(\log x)^{z-1}(\log X)^7}{n(\log x)^2} \ll \frac{(\log X)^8}{(\log x)^{3-z}} \ll \frac{1}{(\log x)^{2-\varepsilon}}.$$

Thus, all error terms from Theorems 2.2 and 2.3, as well as the term $-zy/\log y$, contribute to the right-hand side of (7.4) an amount that can be absorbed by the error term of (7.3). We denote by $S_j^*$ $(j=1,2)$ the contribution of the associated main terms.

Hence, writing $u_n := (\log x/n)/\log nt = v_n - 1$, we have

$$S_1^* := \sum_{\substack{n\in\mathcal{D}(t)\\ n^2t\leqslant X}} \frac{z^{\nu(n)}(\log x/n)^{z-1}}{n}\Big\{h_\nu(nt,z)\Big(1+\frac{(z-1)J_\nu(nt,z)}{\log x/n}\Big)\Big\}$$

$$S_2^* := \sum_{\substack{n\in\mathcal{D}(t)\\ X<n^2t\leqslant x}} \frac{z^{\nu(n)}}{n\log nt}\Big\{\omega_z(u_n) - \frac{(z-1)\alpha(z)(u_n)^{z-2}}{\log nt}\Big\}.$$

At this stage, we extend the summation domains by adding quantities that contribute $\ll (\log t)/(\log x)^{2-\varepsilon}$. To $S_1^*$, we add

$$R_1^* := \sum_{\substack{n\in\mathcal{D}(t)\\ X<n^2t\leqslant x}} \frac{z^{\nu(n)}(\log x/n)^{z-1}}{n}\Big\{h_\nu(nt,z) - \frac{\alpha(z)}{(\log nt)^z}\Big\} \ll \frac{1}{(\log x)^2},$$

by (2.13).

Let $\beta(z) := z\mathrm{e}^{-\gamma z}/\Gamma(z-1) = z(z-1)\alpha(z)$. To $S_2^*$, we add

$$R_2^* := \sum_{\substack{n\in\mathcal{D}(t)\\ n^2t\leqslant X}} \frac{z^{\nu(n)}}{n\log nt}\big\{\omega_z(u_n) - \alpha(z)u_n^{z-1} - \beta(z)u_n^{z-2}\big\} \ll \frac{\log t}{(\log x)^{2-\varepsilon}},$$

by (2.5).

Using the identity $(\log x/n)^{z-1}/(\log nt)^z = u_n^{z-1}/\log nt$, we get, by gathering our estimates,

$$\frac{1}{x}\sum_{m\leqslant x} z^{\nu(m)} = \Delta(x,t,z) + A + B - C + O\Big(\frac{\log t}{(\log x)^{2-\varepsilon}}\Big), \tag{7.5}$$

with

$$A := \sum_{\substack{n\in\mathcal{D}(t)\\ n^2t\leqslant x}} \frac{z^{\nu(n)}}{n}\Big\{\Big(\log\frac{x}{n}\Big)^{z-1} h_\nu(nt,z) + \frac{\omega_z(u_n) - \alpha(z)u_n^{z-1}}{\log nt}\Big\},$$

$$B := \sum_{\substack{n\in\mathcal{D}(t)\\ n^2t\leqslant X}} \frac{z^{\nu(n)}}{n}\Big\{\Big(\log\frac{x}{n}\Big)^{z-2}(z-1)h_\nu(nt,z)J_\nu(nt,z) - \frac{\beta(z)u_n^{z-2}}{\log nt}\Big\},$$

$$C := (z-1)\alpha(z)\sum_{\substack{n\in\mathcal{D}(t)\\ X<n^2t\leqslant x}} \frac{z^{\nu(n)}u_n^{z-2}}{n(\log nt)^2}.$$

In the first term inside curly brackets of $A$, we now replace $\log(x/n)$ by $\log xt$, and in parallel replace $u_n$ by $v_n$ in the third term. The error thus introduced is

$$\begin{aligned}
&\sum_{\substack{n\in\mathcal{D}(t)\\ n^2t\leqslant x}} \frac{z^{\nu(n)}}{n}\Big\{h_\nu(nt,z) - \frac{\alpha(z)}{(\log nt)^z}\Big\}\Big\{(\log xt)^{z-1} - (\log x/n)^{z-1}\Big\}\\
&\qquad = (\log xt)^{z-1}\sum_{\substack{n\in\mathcal{D}(t)\\ n^2t\leqslant x}} \frac{z^{\nu(n)}}{n}\Big\{h_\nu(nt,z) - \frac{\alpha(z)}{(\log nt)^z}\Big\}\Big\{1 - \Big(1 - \frac{1}{v_n}\Big)^{z-1}\Big\}\\
&\qquad = a_{z,t}(\log xt)^{z-2} + O\big((\log x)^{z-3}\big)
\end{aligned}$$

with

$$a_{z,t} := (z-1)\sum_{n\in\mathcal{D}(t)} \frac{z^{\nu(n)}\log nt}{n}\Big\{h_\nu(nt,z) - \frac{\alpha(z)}{(\log nt)^z}\Big\},$$

which is a convergent series by Lemma 2.5.

Next, we extend the summation in $A$ thus modified, to all $n\in\mathcal{D}(t)$. Since $\omega_z(u_n) = 0$ when $n^2t > x$, this introduces an error $\ll 1/(\log x)^2$, by (2.13). Since, by (3.2),

$$\sum_{n\in\mathcal{D}(t)} \frac{z^{\nu(n)}}{n}(\log xt)^{z-1}h_\nu(nt,z) = \lambda_{\nu,0}(z)(\log xt)^{z-1},$$

we see that in (7.5) we can replace $A$ by

$$A^* := \lambda_{\nu,0}(z)(\log xt)^{z-1} - a_{z,t}(\log xt)^{z-2} - \delta(x,t,z).$$

Estimating the left-hand side of (7.5) by Lemma 2.4, we get

$$(\log x)^{z-1}\Big\{\lambda_{\nu,0}(z) + \frac{\lambda_{\nu,1}(z)}{\log x}\Big\} = \Delta(x,t,z) + A^* + B - C + O\Big(\frac{\log t}{(\log x)^{2-\varepsilon}}\Big).$$

Moreover, as defined in Lemma 2.4, $\lambda_{\nu,1}(z)\ll z-1$.

To complete the proof of (7.3), it therefore only remains to show that

$$B - C = b_{z,t}(\log xt)^{z-2} + O\Big(\frac{\log t}{(\log x)^{2-\varepsilon}}\Big),$$

with $b_{z,t} \ll z-1$. To this end, we write $B-C=F-G$, where

$$\frac{F}{z-1} := \sum_{\substack{n\in\mathcal{D}(t)\\ n^2t\leqslant X}} \frac{z^{\nu(n)}}{n}\Big(\log\frac{x}{n}\Big)^{z-2}\Big\{h_\nu(nt,z)(J_\nu(nt,z)+1)-\frac{z\alpha(z)}{(\log nt)^{z-1}}\Big\}$$

and

$$\begin{aligned}\frac{G}{z-1} &:= \sum_{\substack{n\in\mathcal{D}(t)\\ n^2t\leqslant X}} \frac{z^{\nu(n)}h_\nu(nt,z)}{n}\Big(\log\frac{x}{n}\Big)^{z-2} + \sum_{\substack{n\in\mathcal{D}(t)\\ X<n^2t\leqslant x}} \frac{\alpha(z)z^{\nu(n)}}{n(\log nt)^z}\Big(\log\frac{x}{n}\Big)^{z-2}\\ &= \sum_{\substack{n\in\mathcal{D}(t)\\ n^2t\leqslant x}} \frac{z^{\nu(n)}h_\nu(nt,z)}{n}\Big(\log\frac{x}{n}\Big)^{z-2} + O\Big(\frac{1}{(\log x)^2}\Big),\end{aligned}$$

in view of (2.13).

In these two remaining sums, we first replace $\log x/n$ by $\log xt$, which alters $F$ by $\ll (\log x)^{z-3}$ and $G$ by $\ll (\log t)(\log_2 x)/(\log x)^2$. Then we extend both summations from $n^2t\leqslant x$ to all $n\in\mathcal{D}(t)$. We find that

$$F = f_{z,t}(\log xt)^{z-2} + O\Big(\frac{1}{\log x)^2}\Big),\quad G = g_{z,t}(\log xt)^{z-2} + O\Big(\frac{\log t}{(\log x)^{2-\varepsilon}}\Big),$$

where $f_{z,t}$, $g_{z,t} \ll z-1$. This completes the proof. □

Our next lemma provides a functional equation for $x\mapsto\Delta(x,t,z)$. We recall notation (7.1), (7.2), and define

$$\mathfrak{d}_{z,t} := \alpha(z)\int_1^\infty \frac{\Delta(y,t,z)}{y(\log yt)^z}\,\mathrm{d}y \qquad (t\geqslant 2,\ |z|=1,\ \operatorname{Re} z\geqslant 1/2), \tag{7.6}$$

the integral being convergent since $|\Delta(y,t,z)|\leqslant 1$ in the prescribed domain.

**Lemma 7.2.** *Let $\varepsilon>0$. For $z\in\mathbb{C}$, $|z|=1$, $\operatorname{Re} z\geqslant 1/2$, $x\geqslant 1$, $t\geqslant 2$, we have*

$$\begin{aligned}\Delta(x,t,z) &= \mathfrak{d}_{z,t}(\log x)^{z-1} + b_{z,t}(\log xt)^{z-2}\\ &\quad - \int_1^\infty \frac{\Delta(y,t,z)}{y\log yt}\,\omega_z\Big(\frac{\log xt}{\log yt}-1\Big)\,\mathrm{d}y + O\Big(\frac{\log t}{(\log x)^{1-\varepsilon}\log xt}\Big),\end{aligned} \tag{7.7}$$

*where $b_{z,t}\ll (z-1)(\log t)^{1-z}$. Moreover*

$$\mathfrak{d}_{z,t} = \lambda_{\nu,0}(z) + O(1/(\log t)^z). \tag{7.8}$$

*Proof.* Let us first assume $x\geqslant t\geqslant 2$. Consider (7.3). By Abel summation, we have

$$\begin{aligned}\delta(x,t,z) &= \int_1^\infty \Delta(y,t,z)\left\{\alpha(z)v_y^{z-1}-\omega_z(v_y-1)\right\}\frac{1+1/\log yt}{y\log yt}\,\mathrm{d}y\\ &\quad+\int_1^\infty \Delta(y,t,z)\frac{v_y\left\{\alpha(z)(z-1)v_y^{z-2}-\omega_z'(v_y-1)\right\}}{y(\log yt)^2}\,\mathrm{d}y\\ &\quad+O\left(\frac{\log t}{(\log xt)^2}\right),\end{aligned} \tag{7.9}$$

where the error term arises from the discontinuity of $\omega_z(u)$ at $u=1$.

We claim that the contribution from the term $1/\log yt$ in the first integral, as well as the entire second integral in (7.9), can be estimated as

$$e_{z,t}(\log xt)^{z-2} + O\Big(\frac{\log t}{(\log x)^{2-\varepsilon}}\Big)$$

with $e_{z,t} \ll (z-1)(\log t)^{1-z}$. To this end, we shall make frequent use of the estimate

$$|\Delta(y,t,z)| \leqslant \frac{D(y,t)}{y} \ll \frac{\log t}{\log yt} \qquad (y \geqslant 1,\ t \geqslant 2),$$

as well as of the estimates in Lemma 2.1.

The contribution of the term $1/\log yt$ to the range $(x,\infty)$ of the first integral of (7.9) equals

$$\int_x^\infty \frac{\alpha(z)\Delta(y,t,z)}{y(\log yt)^2}\Big(\frac{\log xt}{\log yt}\Big)^{z-1}\,\mathrm{d}y \ll \int_x^\infty \frac{(\log t)(\log xt)^{z-1}}{y(\log yt)^{z+2}}\,\mathrm{d}y \ll \frac{\log t}{(\log xt)^2}.$$

By (2.5), the contribution to the range $[1,x]$ is

$$\int_1^x \frac{\Delta(y,t,z)}{y(\log yt)^2}\big\{-\alpha_1(z)v_y^{z-2}+O(v_y^{z-3})\big\}\,\mathrm{d}y,$$

where $\alpha_1(z) := \mathrm{e}^{-\gamma z}(z-1)/\Gamma(z-1)$. Now, the main term of the last integral equals

$$-\alpha_1(z)\int_1^\infty \frac{\Delta(y,t,z)v_y^{z-2}}{y(\log yt)^2}\,\mathrm{d}y + O\Big(\frac{\log t}{(\log xt)^2}\Big) =: e^*_{z,t}(\log xt)^{z-2} + O\Big(\frac{\log t}{(\log xt)^2}\Big),$$

say, while the error term contributes

$$\ll \int_1^x \frac{(\log t)(\log xt)^{z-3}}{y(\log yt)^z}\,\mathrm{d}y \ll \int_1^x \frac{\log t}{(\log xt)^2 y\log yt}\,\mathrm{d}y \ll \frac{\log t}{(\log xt)^{2-\varepsilon}}.$$

The second integral in (7.9) can be estimated in much the same way, with the help of (2.6) instead of (2.5).

The contribution from the term $\alpha(z)v_y^{z-1}$ to the remaining part of the first integral in (7.9) is

$$\alpha(z)(\log xt)^{z-1}\int_1^\infty \frac{\Delta(y,t,z)}{y(\log yt)^z}\,\mathrm{d}y = \mathfrak{d}_{z,t}(\log xt)^{z-1}, \tag{7.10}$$

with notation (7.6).

We now prove (7.8). By Lemmas 3.3 and 2.5, we have

$$\lambda_{\nu,0}(z) = \alpha(z)\sum_{n\in\mathcal{D}(t)} \frac{z^{\nu(n)}}{n(\log nt)^z} + O\Big(\frac{1}{\log t}\Big),$$

and, by Abel summation, we get that the last sum is

$$\mathfrak{d}_{z,t} + z\alpha(z)\int_1^\infty \frac{\Delta(y,t,z)}{y(\log yt)^{z+1}}\,\mathrm{d}y = \mathfrak{d}_{z,t} + O\Big(\frac{1}{(\log t)^z}\Big).$$

This establishes (7.8).

The above computations open the possibility of combining (7.10) with the second term on the right-hand side of (7.3):

$$\begin{aligned}
&\mathfrak{d}_{z,t}(\log xt)^{z-1} + \lambda_{\nu,0}(z)\{(\log x)^{z-1}-(\log xt)^{z-1}\}\\
&\quad = \{\mathfrak{d}_{z,t}-\lambda_{\nu,0}(z)\}\{(\log xt)^{z-1}-(\log x)^{z-1}\} + \mathfrak{d}_{z,t}(\log x)^{z-1}\\
&\quad = \{\mathfrak{d}_{z,t}-\lambda_{\nu,0}(z)\}(z-1)(\log xt)^{z-2}\log t + \mathfrak{d}_{z,t}(\log x)^{z-1} + O\Big(\frac{\log t}{(\log x)^2}\Big)\\
&\quad = b^*_{z,t}(\log xt)^{z-2} + \mathfrak{d}_{z,t}(\log x)^{z-1} + O\Big(\frac{\log t}{(\log x)^2}\Big),
\end{aligned}$$

where $b^*_{z,t} = \{\mathfrak{d}_{z,t}-\lambda_{\nu,0}(z)\}(z-1)\log t \ll (z-1)(\log t)^{1-z}$. This proves the result in the case $x \geqslant t \geqslant 2$, with $b_{z,t} = c_{z,t}+e_{z,t}+b^*_{z,t} \ll (z-1)(\log t)^{1-z}$.

If $1 \leqslant x \leqslant t$, the integral in Lemma 7.2 vanishes, the term $b_{z,t}(\log xt)^{z-2}$ is absorbed by the error term, and the result follows from Lemma 2.4 and (7.8). □

**Proposition 7.3.** *Let $0 < \kappa < 1$. There exists $\delta > 0$, such that the following statement holds. Assume that for $x \geqslant 1$, $|z| = 1$ and $|z-1| \leqslant \delta$, $\Delta(x,t,z)$ satisfies $|\Delta(x,t,z)| \leqslant 1$ and*

$$\begin{aligned}\Delta(x,t,z) &= \mathfrak{d}_{z,t}(\log x)^{z-1} + b_{z,t}(\log xt)^{z-2}\\ &\quad - \int_1^\infty \frac{\Delta(y,t,z)}{y\log yt}\omega_z\left(\frac{\log xt}{\log yt}-1\right)\mathrm{d}y + O\left(\frac{\log t}{(\log x)^\kappa \log xt}\right),\end{aligned} \tag{7.11}$$

*where $b_{z,t} \ll (\log t)^{1-\kappa}$. Then, for $x \geqslant 1$, $|z| = 1$ and $|z-1| \leqslant \delta$, we have*

$$\Delta(x,t,z) = \alpha_{z,t}d_z(v) + \frac{\beta_{z,t}}{v} + O\left(\frac{\log t}{(\log x)^\kappa \log xt}\right),$$

*where $v = (\log x)/\log t$, $\alpha_{z,t} = \mathfrak{d}_{z,t}(\log t)^{z-1} + O\big(1/(\log t)^\kappa\big)$ and $\beta_{z,t} \ll 1/(\log t)^\kappa$.*

*Proof.* The parameters $t$ and $z$ being fixed in their prescribed domains, define $w$, $u$, and $G$ by

$$\mathrm{e}^w = \frac{\log xt}{\log t}, \qquad \mathrm{e}^u = \frac{\log yt}{\log t}, \qquad G(w) = \Delta(x,t,z).$$

The assumption on $\Delta(x,t,z)$ is equivalent to

$$G(w) = \mathfrak{a}_{z,t}(\mathrm{e}^w-1)^{z-1} + \mathfrak{b}_{z,t}\mathrm{e}^{w(z-2)} - \int_0^w G(u)M_z(w-u)\,\mathrm{d}u + \mathfrak{E}(w),$$

where

$$\mathfrak{a}_{z,t} = \mathfrak{d}_{z,t}(\log t)^{z-1},\ \mathfrak{b}_{z,t} = b_{z,t}(\log t)^{z-2},\ \mathfrak{E}(w) \ll E_t(w) := \frac{1}{\mathrm{e}^w(e^w-1)^\kappa(\log t)^\kappa}\cdot$$

Now $|G(w)| \leqslant 1$ and $|M_z(w)| \leqslant 2(\mathrm{e}^u-1)$ for $|z| \leqslant 2$, by (2.7). Multiplying by $\mathrm{e}^{-sw}$, where $\operatorname{Re} s > 1$, and integrating with respect to $w \geqslant 0$, leads to the equation of Laplace transforms

$$\widehat{G}(s) = \mathfrak{a}_{z,t}\frac{\Gamma(s+1-z)\,\Gamma(z)}{\Gamma(s+1)} + \frac{\mathfrak{b}_{z,t}}{s+2-z} - \widehat{G}(s)\widehat{M_z}(s) + \widehat{\mathfrak{E}}(s),$$

as in (5.7). Solving for $\widehat{G}(s)$, we find that

$$\widehat{G}(s) = \mathfrak{a}_{z,t}\widehat{q_z}(s) + \mathfrak{b}_{z,t}\widehat{F_z}(s) + \widehat{\mathfrak{E}}(s)\big\{1-\widehat{H_z}(s)\big\}, \tag{7.12}$$

where $q_z(w) = d_z(e^w-1)$, as defined in (5.1), and

$$\widehat{F_z}(s) := \frac{1}{(s+2-z)(1+\widehat{M_z}(s))}, \qquad \widehat{H_z}(s) := 1 - \frac{1}{1+\widehat{M_z}(s)}\cdot \tag{7.13}$$

Note that $\widehat{F_z}(s)$ and $\widehat{H_z}(s)$ are the Laplace transforms of the functions $F_z(w)$ and $H_z(w)$, defined by

$$\begin{aligned}F_z(w) &= \mathrm{e}^{w(z-2)} - \int_0^w F_z(u)M_z(w-u)\,\mathrm{d}u,\\ H_z(w) &= M_z(w) - \int_0^w H_z(u)M_z(w-u)\,\mathrm{d}u,\end{aligned}$$

and initial condition $M_z(w) = 0$ for $0 \leqslant w < \log 2$: indeed these equations uniquely define $F_z(w)$ and $H_z(w)$ on successive intervals $[(k-1)\log 2, k\log 2]$, $k \geqslant 1$. As $|M_z(w)| \leqslant 2(\mathrm{e}^u - 1)$ for $|z| \leqslant 2$, by (2.7), it follows by induction on $k \geqslant 1$, that for $w \in [(k-1)\log 2, k \log 2]$, we have the preliminary estimates $|F_z(w)| \leqslant \mathrm{e}^{2w}$ and $|H_z(w)| \leqslant \mathrm{e}^{2w}$.

Inverting equation (7.12), we obtain

$$G(w) = \mathfrak{a}_{z,t} q_z(w) + \mathfrak{b}_{z,t} F_z(w) + \mathfrak{E}(w) - \int_0^w \mathfrak{E}(u) H_z(w-u)\,\mathrm{d}u. \tag{7.14}$$

Carrying (4.10) back into (7.13), we find that, for $|\tau| \geqslant 1$, $\sigma \geqslant -\frac{5}{2}$, we have

$$\widehat{F_z}(s) = \frac{1}{s+2-z} + O\Big(\frac{1}{\tau^2}\Big), \qquad \widehat{H_z}(s) = \frac{z}{s2^s} + O\Big(\frac{1}{\tau^2}\Big). \tag{7.15}$$

Hence we may write

$$F_z(w) = \frac{1}{2\pi i}\int_{3-i\infty}^{3+i\infty} \widehat{F_z}(s)\mathrm{e}^{ws}\,\mathrm{d}s, \quad H_z(w) = \frac{1}{2\pi i}\int_{3-i\infty}^{3+i\infty} \widehat{H_z}(s)\mathrm{e}^{ws}\,\mathrm{d}s.$$

Lemma 4.5 shows that for $\operatorname{Re} s \geqslant -5/2$, the poles of $\widehat{F_z}(s)$ and $\widehat{H_z}(s)$ are at $s = s_0(z)$, $-1$, $-2$, since, according to Corollary 4.2, the zero of the factor $(s+2-z)$ in (7.13) is canceled by a pole of $\widehat{M_z}(s)$. With the estimates (7.15), we can truncate the integrals at height $|\operatorname{Im} s| = \mathrm{e}^{5w}$ and use the residue theorem to move the integration line to the left as far as $\operatorname{Re} s = -5/2$. After some standard calculations, this shows that

$$F_z(w) = r_{0,z}\mathrm{e}^{s_0(z)w} + r_{1,z}\mathrm{e}^{-w} + O(\mathrm{e}^{-2w}), \tag{7.16}$$

$$H_z(w) = r^*_{0,z}\mathrm{e}^{s_0(z)w} + r^*_{1,z}\mathrm{e}^{-w} + O(\mathrm{e}^{-2w}), \tag{7.17}$$

where the $r$'s are the residues at the corresponding poles. With the help of (7.17), we can approximate the integral in (7.14) as

$$\begin{aligned}\int_0^w \mathfrak{E}(u)H_z(w-u)\,\mathrm{d}u &= r^*_{0,z}\mathrm{e}^{s_0(z)w}\int_0^w \mathfrak{E}(u)\mathrm{e}^{-s_0(z)u}\,\mathrm{d}u + r^*_{1,z}\mathrm{e}^{-w}\int_0^w \mathfrak{E}(u)\mathrm{e}^{u}\,\mathrm{d}u \\ &\quad + O\Big(\int_0^w E_t(u)\mathrm{e}^{-2(w-u)}\,\mathrm{d}u\Big) \\ &= r^*_{0,z}\mathrm{e}^{s_0(z)w}\widehat{\mathfrak{E}}(s_0(z)) + r^*_{1,z}\mathrm{e}^{-w}\widehat{\mathfrak{E}}(-1) + O\big(E_t(w)\big).\end{aligned}$$

Inserting this estimate and (7.16) into (7.14), we conclude that

$$G(w) = \mathfrak{a}_{z,t}q_z(w) + \xi_{z,t}\mathrm{e}^{s_0(z)w} + \beta_{z,t}\mathrm{e}^{-w} + O\big(|\mathfrak{b}_{z,t}|\mathrm{e}^{-2w} + E_t(w)\big),$$

where

$$\xi_{z,t} := \mathfrak{b}_{z,t}r_{0,z} - r^*_{0,z}\widehat{\mathfrak{E}}(s_0(z)), \tag{7.18}$$

$$\beta_{z,t}; = \mathfrak{b}_{z,t}r_{1,z} - r^*_{1,z}\widehat{\mathfrak{E}}(-1). \tag{7.19}$$

Now $\mathfrak{b}_{z,t} = b_{z,t}(\log t)^{z-2} \ll 1/(\log t)^\kappa$, while

$$\widehat{\mathfrak{E}}(s_0(z)) \ll 1/(\log t)^\kappa, \quad \widehat{\mathfrak{E}}(-1) \ll 1/(\log t)^\kappa,$$

since $\mathfrak{E}(w) \ll E_t(w)$. It follows that $\xi_{z,t} \ll 1/(\log t)^\kappa$ and $\beta_{z,t} \ll 1/(\log t)^\kappa$. Together with Theorem 5.1, this yields

$$G(w) = \alpha_{z,t}q_z(w) + \beta_{z,t}\mathrm{e}^{-w} + O\big(E_t(w)\big),$$

where, recalling (5.4),

$$\alpha_{z,t} = \mathfrak{a}_{z,t} + \xi_{z,t}/C_z = \mathfrak{d}_{z,t}(\log t)^{z-1} + O(1/(\log t)^\kappa). \tag{7.20}$$

□

**Lemma 7.4.** *For* $s, z \in \mathbb{C}$, $\operatorname{Re} s > 1$, $|z-1| \leqslant 1/2$, *we have*

$$1 = \sum_{n\in\mathcal{B}} \frac{z^{\Omega(n)}}{n^s} \prod_{p\leqslant\vartheta(n)} \Big(1-\frac{z}{p^s}\Big), \quad 1 = \sum_{n\in\mathcal{B}} \frac{z^{\omega(n)}}{n^s} \prod_{p\leqslant\vartheta(n)} \Big(1+\frac{z}{p^s-1}\Big)^{-1}. \tag{7.21}$$

*Proof.* The canonical representation made explicit in the proof of Lemma 3.1 leads to the identity

$$\begin{aligned}\prod_{p\geqslant2}\Big(1-\frac{z}{p^s}\Big)^{-1} = \sum_{m\geqslant1}\frac{z^{\Omega(m)}}{m^s} &= \sum_{n\in\mathcal{B}}\frac{z^{\Omega(n)}}{n^s}\sum_{P^-(r)>\vartheta(n)}\frac{z^{\Omega(r)}}{r^s}\\ &= \sum_{n\in\mathcal{B}}\frac{z^{\Omega(n)}}{n^s}\prod_{p>\vartheta(n)}\Big(1-\frac{z}{p^s}\Big)^{-1},\end{aligned}$$

for $\operatorname{Re} s > 1$. The result now follows from dividing by the product on the left-hand side. The identity involving $\omega(n)$ is derived similarly. □

For $x \geqslant 1$, $z \in \mathbb{C}$, let us introduce the notation

$$B(x,z) := \sum_{n\in\mathcal{B}(x)} z^{\nu(n)}. \tag{7.22}$$

**Lemma 7.5.** *Let* $a < 1$. *and assume* $\max(2,n) \leqslant \vartheta(n) \ll n\exp\{(\log n)^a\}$ *for all* $n \geqslant 1$. *For suitable* $\delta > 0$ *and* $|z-1| \leqslant \delta$, $|z| = 1$, *the following implication holds. If* $c_z$, $w \in \mathbb{C}$, $\varepsilon > 0$, *are such that*

$$B(x,z) = c_z x(\log x)^w \Big(1 + O\Big(\frac{1}{(\log x)^\varepsilon}\Big)\Big\} \qquad (x \geqslant 2), \tag{7.23}$$

*then* $\operatorname{Re} w \leqslant \operatorname{Re} z - 2$.

*Proof.* We only consider $\nu = \Omega$, the case $\nu = \omega$ being very similar.

Let us argue by contradiction and assume $\operatorname{Re} w > \operatorname{Re} z - 2$ in (7.23). Differentiating (7.21) with respect to $s$, we have

$$\sum_{n\in\mathcal{B}} \frac{z^{\Omega(n)}}{n^s}\Bigg(\sum_{p\leqslant\vartheta(n)}\frac{z\log p}{p^s-z} - \log n\Bigg)\prod_{p\leqslant\vartheta(n)}\Big(1-\frac{z}{p^s}\Big) = 0. \tag{7.24}$$

Let $F_{N,z}(s)$ (resp. $G_{N,z}(s)$) denote the contribution to the right-hand side of (7.24) from $n \leqslant N$ (resp. $n > N$). The following calculations are very similar to those in [18, Lemmas 3 and 4], where more details are given. With $s = 1 + 1/(\log N)^2$, we find that, as $N \to \infty$, with notation (2.10),

$$\begin{aligned}F_{N,z}(s) &\sim \sum_{n\in\mathcal{B}(N)}\frac{z^{\Omega(n)}}{n}\Bigg(\sum_{p\leqslant\vartheta(n)}\frac{z\log p}{p-z} - \log n\Bigg)\prod_{p\leqslant\vartheta(n)}\Big(1-\frac{z}{p}\Big)\\ &\sim \sum_{n\in\mathcal{B}(N)}\frac{z^{\Omega(n)}\mathrm{e}^{-\gamma z}(z-1)}{H_\Omega(1,z)n(\log n)^{z-1}} \sim \frac{\mathrm{e}^{-\gamma z}(z-1)}{H_\Omega(1,z)}\int_2^N \frac{B(y,z)}{y^2(\log y)^{z-1}}\,\mathrm{d}y\\ &\sim \frac{c_z(z-1)\mathrm{e}^{-\gamma z}(\log N)^{w+2-z}}{H_\Omega(1,z)(w+2-z)},\end{aligned}$$

since $\operatorname{Re}(w-z+2)>0$. On the other hand, for the same value of $s$,

$$G_{N,z}(s)=\sum_{\substack{n\in\mathcal{B}\\ n>N}}\frac{z^{\Omega(n)}}{n^s}\Bigg(\sum_{p\leqslant\vartheta(n)}\frac{z\log p}{p^s-z}-\log n\Bigg)\prod_{p\leqslant\vartheta(n)}\Big(1-\frac{z}{p^s}\Big)$$
$$\sim\sum_{\substack{n\in\mathcal{B}\\ n>N}}\frac{z^{\Omega(n)}\mathrm{e}^{-\gamma z+zT((s-1)\log n)}}{H_\Omega(s,z)n^s(\log n)^z}\Big(z\frac{1-n^{1-s}}{s-1}-\log n\Big),$$

where $T(s)$ is defined in (4.5). Abel summation, followed by the change of variables $u=(s-1)\log n$ leads to

$$G_{N,z}(s)\sim\frac{c_z(s-1)^{z-2-w}}{\mathrm{e}^{\gamma z}H_\Omega(1,z)}\int_{1/\log N}^\infty u^{w-z+1}\left\{zT'(u)-1\right\}\mathrm{e}^{zT(u)-u}\,\mathrm{d}u.$$

If $\operatorname{Re}(w-z+1)>-1$, the last integral converges to some value $I_{w,z}$, say. Since $s-1=1/(\log N)^2$, $F_{N,z}(s)$ and $G_{N,z}(s)$ do not have the same order of magnitude if $I_{w,z}\neq 0$. As $F_{N,z}(s)+G_{N,z}(s)=0$ by (7.24), we must have $I_{w,z}=0$. However, integration by parts taking account of (4.1) imply that, for $\operatorname{Re}(w-z+1)>0$, we have

$$I_{w,z}=-\mathrm{e}^{\gamma z}(w-z+1)f_z(w),$$

where $f_z$ is defined in (4.3). This formula extends to all $w\in\mathbb{C}$ by meromorphic continuation. Since $f_z(w)$ has a pole at $w=z-1$ and $f_z(w)\neq 0$ if $\operatorname{Re}w>\operatorname{Re}z-2$, by Lemma 4.4, this contradicts the vanishing of $I_{w,z}$. Therefore, we must have $\operatorname{Re}w\leqslant\operatorname{Re}z-2$. □

## 8. Proof of Theorem 6.1: completion of the argument

Lemma 7.2 shows that the assumptions of Proposition 7.3 are fulfilled, with $\kappa=1-\varepsilon/2$. We may hence deduce from Lemma 7.5 that $\beta_{z,t}$ must vanish, because $\operatorname{Re}s_0(z)<-1$, by Lemma 4.4. This implies estimate (6.2) with

$$\eta_{z,t}=\alpha_{z,t}(\log t)^{1-z}=\mathfrak{a}_{z,t}(\log t)^{1-z}+\frac{\xi_{z,t}(\log t)^{1-z}}{C_z}=\mathfrak{d}_{z,t}+O\Big(\frac{1}{(\log t)^{1-\varepsilon}}\Big),\tag{8.1}$$

by (7.20), where $C_z$ is defined in (5.4) and $\xi_{z,t}$ in (7.18).

It remains to show that $\eta_{z,t}=\eta_t+O(z-1)$. Combining (6.4) with (5.2) and (5.3), we get for $z=\mathrm{e}^{i\varphi}$, $|z-1|\leqslant\delta$, $\varphi^2\mu_{x,t}\ll 1$,

$$\Delta(x,t,z)=\eta_tC\mathrm{e}^{-\log_2(xt)+\log_2 t+(z-1)\mu_{x,t}}\Big\{1+O\Big(\frac{1}{\log x}+\frac{1}{v^2}+|\varphi|+\varphi^2\mu_{x,t}^2\Big)\Big\},$$

while (6.2) and Theorem 5.1 yield

$$\Delta(x,t,z)=\eta_{z,t}C_z\mathrm{e}^{(z-1)\log_2 t+s_0(z)(\log_2(xt)-\log_2 t)}\Big\{1+O\Big(\frac{1}{\log x)^{1-\varepsilon}}+\frac{1}{v^2}\Big)\Big\}.$$

We estimate $s_0(z)$ by Lemma 4.4 and insert the estimate $C_z=C+O(z-1)$ from (5.5) to obtain

$$\eta_{z,t}=\eta_t\Big\{1+O\Big(|\varphi|+\varphi^2\mu_{x,t}^2+\frac{1}{v^2}+\frac{1}{(\log x)^{1-\varepsilon}}\Big)\Big\}.$$

The desired estimate (6.3) follows with $\log x=\mathrm{e}^{2/\sqrt{|\varphi|}}$, provided $\log t\leqslant\sqrt{\log x}$. If $\log t>\sqrt{\log x}=\exp\big(1/\sqrt{|\varphi|}\big)$, then (7.8) and (8.1) yield

$$\eta_{z,t}=\mathfrak{d}_{z,t}+O\Big(\frac{1}{(\log t)^{1-\varepsilon}}\Big)=\lambda_{\nu,0}(z)+O(\varphi)=1+O(\varphi)=\eta_t\{1+O(\varphi)\}.$$

## 9. Proof of Theorem 1.2

Let $F_x(y)$ denote the left-hand side of (1.1) and let

$$\psi_x(\tau) := \frac{1}{D(x,t)} \sum_{n\in\mathcal{D}(x,t)} \mathrm{e}^{i\tau(\nu(n)-\mu_{x,t})/\sigma_{x,t}} \qquad (\tau \in \mathbb{R}).$$

With $z = \mathrm{e}^{i\varphi}$, we have $z-1 = i\varphi - \frac{1}{2}\varphi^2 + O\big(\varphi^3\big)$ and $(z-1)^2 = -\varphi^2 + O\big(\varphi^3\big)$. Lemma 4.4 and Corollary 6.3 yield, uniformly for $x \geqslant t \geqslant 2$, $|\varphi| \leqslant \delta$,

$$\begin{aligned} D(x,t,\mathrm{e}^{i\varphi})/D(x,t) = \{1+O(\varphi)\}&\frac{(\log x)^{Ci\varphi - V\varphi^2/2+O(\varphi^3)}}{(\log t)^{(C-1)i\varphi+(1-V)\varphi^2/2+O(\varphi^3)}} \\ &+ O\Big(\frac{1}{(\log x)^{1-\varepsilon}}\Big). \end{aligned} \tag{9.1}$$

With

$$T := \delta\, \sigma_{x,t}, \quad \varphi := \tau/\sigma_{x,t},$$

we obtain

$$\psi_x(\tau) \ll (\log x)^{-V\varphi^2/3}(\log t)^{-(1-V)\varphi^2/3} = \mathrm{e}^{-\tau^2/3} \qquad (|\tau| \leqslant T), \tag{9.2}$$

provided $\delta$ is sufficiently small. We will use this estimate for $T^{1/3} < |\tau| \leqslant T$.

When $|\tau| \leqslant T^{1/3}$, (9.1) yields

$$\psi_x(\tau) = \mathrm{e}^{-\tau^2/2}\left\{1 + O\left(\frac{|\tau|+|\tau|^3}{T}\right)\right\} + O\Big(\frac{1}{(\log x)^{1-\varepsilon}}\Big) \quad (|\tau| \leqslant T^{1/3}), \tag{9.3}$$

which we shall use in the range $1/T^2 < |\tau| \leqslant T^{1/3}$.

When $|\tau| \leqslant 1/T^2$, we insert the trivial estimates

$$\mathrm{e}^{iy} = 1 + O(y), \quad |\nu(n) - \mu_{x,t}| \leqslant \nu(n) + \mu_{x,t}$$

in the definition of $\psi_x(\tau)$ to obtain

$$\psi_x(\tau) - 1 \ll \frac{|\tau|}{TD(x,t)} \sum_{n\in\mathcal{D}(x,t)} (\nu(n) + \mu_{x,t}) \ll |\tau| T, \tag{9.4}$$

by Lemma 6.4.

The Berry-Esseen inequality (see e.g. [13, Thm II.7.16]) states that

$$\sup_{y\in\mathbb{R}} |F_x(y) - \Phi(y)| \ll \frac{1}{T} + \int_{-T}^{T} \Big|\psi_x(\tau) - \mathrm{e}^{-\tau^2/2}\Big| \frac{\mathrm{d}\tau}{|\tau|}.$$

We split the integral into three parts, $I_1$, $I_2$ and $I_3$, according to the three different ranges for $|\tau|$ mentioned above. By (9.2), (9.3) and (9.4), we have

$$\begin{aligned} I_1 &\ll \int_{T^{1/3}}^{T} \mathrm{e}^{-\tau^2/3}\frac{\mathrm{d}\tau}{\tau} \ll \frac{1}{T}, \\ I_2 &\ll \int_{1/T^2}^{T^{1/3}} \frac{1+\tau^2}{T}\,\mathrm{e}^{-\tau^2/2}\,\mathrm{d}\tau + \frac{1}{(\log x)^{1-\varepsilon}} \int_{1/T^2}^{T^{1/3}} \frac{\mathrm{d}\tau}{\tau} \ll \frac{1}{T}, \\ I_3 &\ll T\int_{-1/T^2}^{1/T^2} \mathrm{d}\tau \ll \frac{1}{T}. \end{aligned}$$

This completes the proof of Theorem 1.2.

## 10. Summing $z^{\nu(n)}$ over $\mathcal{B}(x)$

Recall notation (7.22).

**Theorem 10.1.** *Let* $0 < a < 1$ *and assume* $\max(2, n) \leqslant \vartheta(n) \ll n \exp\{(\log n)^a\}$ *for all* $n \geqslant 1$*. For suitable* $\delta > 0$ *and uniformly for* $|z| = 1$*,* $|z - 1| < \delta$ *and* $x \geqslant 2$*, we have*

$$B(x, z) = \alpha_z x(\log x)^{s_0(z)} + O\Big(\frac{x}{(\log x)^{2-a}}\Big),$$

*where* $\alpha_z = c_\vartheta + O(z - 1)$ *and* $c_\vartheta$ *is defined in Lemma 3.2.*

From Theorem 10.1 and Lemma 3.2, we immediately derive the following statement.

**Corollary 10.2.** *Let* $0 < a < 1$ *and assume* $\max(2, n) \leqslant \vartheta(n) \ll n \exp\{(\log n)^a\}$ *for all* $n \geqslant 1$*. For suitable* $\delta > 0$ *and uniformly for* $|z| = 1$*,* $|z - 1| < \delta$ *and* $x \geqslant 2$*, we have*

$$B(x, z)/B(x) = \lambda_z (\log x)^{1+s_0(z)} + O\Big(\frac{1}{(\log x)^{1-a}}\Big),$$

*where* $\lambda_z = \alpha_z / c_\vartheta = 1 + O(z - 1)$*.*

Let us now embark on the proof of Theorem 10.1. Define

$$\beta(x, z) := \frac{B(x, z)}{x}.$$

**Lemma 10.3.** *Let* $0 < a < 1$ *and assume* $\max(2, n) \leqslant \vartheta(n) \ll n \exp\{(\log n)^a\}$ *for all* $n \geqslant 1$*. For* $|z| = 1$*,* $\operatorname{Re} z \geqslant 1/2$ *and* $x \geqslant 3$*, we have*

$$\begin{aligned}\beta(x, z) = \sum_{n \in \mathcal{B}} \frac{z^{\nu(n)}}{n \log \vartheta(n)} &\left\{ \frac{\mathrm{e}^{-\gamma z}}{\Gamma(z)} \left( \frac{\log \mathrm{e}x}{\log \vartheta(n)} \right)^{z-1} - \omega_z \left( \frac{\log x/n}{\log \vartheta(n)} \right) \right\} \\ &+ k_z (\log \mathrm{e}x)^{z-2} + O\Big(\frac{1}{(\log \mathrm{e}x)^{2-a}}\Big),\end{aligned} \tag{10.1}$$

*where* $k_z \ll z - 1$ *depends only on* $z$*.*

*Proof.* The proof is almost the same as that of Lemma 7.1, with $\vartheta(n)$ replacing $nt$. We have $\log \vartheta(n) = \log \mathrm{e}n + O((\log \mathrm{e}n)^a)$. The error term $(\log t)/(\log x)^{2-\varepsilon}$ can be replaced by $1/(\log \mathrm{e}x)^{2-a}$ throughout.

The contribution from the term $-zy/\log y$ from Theorem 2.2 can be estimated as

$$\ll \sum_{\substack{n \in \mathcal{B} \\ n\vartheta(n) \leqslant x}} \frac{\vartheta(n)}{\log \vartheta(n)} \ll \frac{x}{(\log x)^{2-a}},$$

as in the proof of [20, th. 4].

While in the proof of Lemma 7.1 we replaced $\log(x/n)$ by $\log xt$, here we replace it by $\log \mathrm{e}x$.

Note that the term $\lambda_{\nu,0}(z)\{(\log x)^{z-1} - (\log xt)^{z-1}\}$ appearing in the statement of Lemma 7.1 is not needed here, since

$$(\log x)^{z-1} - (\log \mathrm{e}x)^{z-1} = -(z - 1)(\log \mathrm{e}x)^{z-2} + O\big(1/(\log \mathrm{e}x)^{2-a}\big).$$

□

**Lemma 10.4.** *Let* $0 < a < 1$ *and assume* $\max(2, n) \leqslant \vartheta(n) \ll n \exp\{(\log n)^a\}$ *for all* $n \geqslant 1$*. For* $|z| = 1$, $\operatorname{Re} z \geqslant 1/2$ *and* $x \geqslant 1$*, we have*

$$\beta(x, z) = \sum_{n \in \mathcal{B}} \frac{z^{\nu(n)}}{n \log \mathrm{e}n} \left\{ \frac{\mathrm{e}^{-\gamma z}}{\Gamma(z)} \left( \frac{\log \mathrm{e}x}{\log \mathrm{e}n} \right)^{z-1} - \omega_z \left( \frac{\log \mathrm{e}x}{\log \mathrm{e}n} - 1 \right) \right\} + k_z (\log \mathrm{e}x)^{z-2} + O((\log \mathrm{e}x)^{a-2}).$$

*Proof.* To derive this result from Lemma 10.3, we will make three changes to the sum on the right-hand side of (10.1). Each of these modifications alters the value of the sum by an expression of type

$$k_z^* (\log \mathrm{e}x)^{z-2} + O((\log \mathrm{e}x)^{a-2}), \tag{10.2}$$

where $k_z^*$ only depends on $z$ and $\vartheta$.

First, we replace the argument of $\omega_z$ by

$$u_n := \frac{\log \mathrm{e}x}{\log \vartheta(n)} - 1,$$

thereby altering the sum by

$$\sum_{n \in \mathcal{B}} \frac{z^{\nu(n)} \{\omega_z(u_n) - \omega_z(u_n + \varepsilon_n)\}}{n \log \vartheta(n)},$$

where

$$\varepsilon_n := 1 - \frac{\log \mathrm{e}n}{\log \vartheta(n)} \ll \frac{1}{(\log \mathrm{e}n)^{1-a}}.$$

Write

$$\omega_z(u_n) - \omega_z(u_n + \varepsilon_n) = -\varepsilon_n \omega_z'(u_n) + O\big(\varepsilon_n^2 \omega_z''(u_n)\big),$$

and use Lemma 2.1 to estimate $\omega_z'(u)$ and $\omega_z''(u)$ when $u_n \geqslant 0$, i.e. $\vartheta(n) \leqslant \mathrm{e}x$. We then extend the sum corresponding to $-\varepsilon_n \omega_z'(u_n)$ from $\vartheta(n) \leqslant \mathrm{e}x$ to all $n \in \mathcal{B}$. Because of the discontinuity of $\omega_z(u)$ at $u = 1$, we also consider the contribution from $n$ where $u_n < 1 \leqslant u_n + \varepsilon_n$ or $u + \varepsilon_n < 1 < u$, which turns out to be $\ll (\log \mathrm{e}x)^{a-2}$. This shows that the first modification has an effect of (10.2).

Second, we replace the first instance of $\log \vartheta(n)$ by $\log \mathrm{e}n$. This alters the sum by

$$\sum_{n \in \mathcal{B}} \frac{\varepsilon_n z^{\nu(n)} \{\alpha(z)(u_n + 1)^{z-1} - \omega_z(u_n)\}}{n \log \mathrm{e}n}.$$

The contribution from $\vartheta(n) > \mathrm{e}x$ is $\ll (\log \mathrm{e}x)^{a-2}$. For $u_n \geqslant 0$, i.e. $\vartheta(n) \leqslant \mathrm{e}x$, we use Lemma 2.1 to estimate $\alpha(z)(u+1)^{z-1} - \omega_z(u)$. We extend the sum corresponding to the main term of this approximation from $\vartheta(n) \leqslant \mathrm{e}x$ to all $n \in \mathcal{B}$, and estimate the contribution from the error term. This shows that the second modification also has an effect of type (10.2).

Our third and last modification consists in replacing $u_n$ by $w_n$, where

$$w_n := \frac{\log \mathrm{e}x}{\log \mathrm{e}n} - 1.$$

Applying the mean value theorem, the implied error is

$$\sum_{n \in \mathcal{B}} \frac{z^{\nu(n)} \{f(w_n) - f(u_n)\}}{n \log \mathrm{e}n} = \sum_{n \in \mathcal{B}} \frac{z^{\nu(n)} f'(u_n^*)(w_n - u_n)}{n \log \mathrm{e}n}, \tag{10.3}$$

where

$$f(u) := \alpha(z)(u + 1)^{z-1} - \omega_z(u),$$

and $u_n^*$ lies between $u_n$ and $w_n$. Estimating each $f$ and $f'$ by Lemma 2.1, we find that, when $u_n \neq w_n$,

$$(u_n^* + 1)^{z-3} = \frac{\kappa_z(n)(\log \mathrm{e}x)^{z-3}}{z-2}\Big\{1 + O\Big(\frac{1}{\log \mathrm{e}x}\Big)\Big),$$

where

$$\kappa_z(n) = \frac{\{\log \vartheta(n)\}^{2-z} - \{\log \mathrm{e}n\}^{2-z}}{1/\log \vartheta(n) - 1/\log \mathrm{e}n}.$$

With this estimate and Lemma 2.1, we conclude that the second sum in (10.3) is also of type (10.2). □

**Lemma 10.5.** *Let* $0 < a < 1$ *and assume* $\max(2, n) \leqslant \vartheta(n) \ll n \exp\{(\log n)^a\}$ *for all* $n \geqslant 1$*. For* $|z| = 1$*,* $\operatorname{Re} z \geqslant 1/2$ *and* $x \geqslant 1$*,*

$$\beta(x,z) = a_z(\log x)^{z-1} + b_z(\log \mathrm{e}x)^{z-2} - \int_1^\infty \frac{\beta(y,z)}{y\log \mathrm{e}y}\,\omega_z\left(\frac{\log \mathrm{e}x}{\log \mathrm{e}y} - 1\right)\mathrm{d}y + O\Big(\frac{x}{(\log \mathrm{e}x)^{2-a}}\Big),$$

*where* $a_z$ *and* $b_z$ *depend only on* $z$*.*

*Proof.* This follows from Lemma 10.4 with Abel summation, in a similar way to the proof of Lemma 7.2. □

**Lemma 10.6.** *Let* $0 < a < 1$ *and assume* $\max(2, n) \leqslant \vartheta(n) \ll n \exp\{(\log n)^a\}$ *for all* $n \geqslant 1$*. Let* $k \in \mathbb{N}$ *be fixed. We have*

$$\sum_{n\in\mathcal{B}(x)} \nu(n)^k \ll_k B(x)(\log_2 x)^k.$$

*Proof.* Put $r_k := (k+2)/\log 2$. The contribution from those $n$ with $\nu(n) \leqslant r_k \log_2 x$ is clearly acceptable. The contribution from those $n$ with $r_k \log_2 x < \nu(n)$ is also acceptable, since $B(x) \asymp x/\log x$ by Lemma 3.2, and, for any fixed $\varepsilon > 0$,

$$|\{n \leqslant x : \Omega(n) \geqslant y\}| \ll_\varepsilon \frac{x\log x}{2^y} \qquad \big(y > (2+\varepsilon)\log_2 x\big),$$

by a theorem of Nicolas [6]—se also [13, th. III.6. & ex. 217]. □

*Proof of Theorem 10.1.* Applying Lemma 10.5 and Proposition 7.3, with $t = \mathrm{e}$, $\kappa = 1 - a$, and taking account of Theorem 5.1, we have

$$\beta(x,z) = \alpha_z(\log x)^{s_0(z)} + \beta_z(\log x)^{-1} + O((\log x)^{a-2}). \tag{10.4}$$

Lemma 7.5 shows that $\beta_z$ vanishes for $|z - 1| < \delta$.

It remains to show that

$$\alpha_z = c_\vartheta + O(z-1). \tag{10.5}$$

Writing $z^{\nu(n)} = e^{i\varphi\nu(n)} = 1 + i\varphi\nu(n) + O\big(\varphi^2\nu(n)^2\big)$, Lemma 10.6 yields

$$B(x,z) = B(x)\big\{1 + i\varphi N_x + O\big(\varphi^2\mu_x^2\big)\big\}, \tag{10.6}$$

where $N_x = B(x)^{-1}\sum_{n\in\mathcal{B}(x)}\nu(n)$ and $\mu_x = C\log_2 x$. We estimate $B(x)$ with Lemma 3.2 and equate the result with (10.4), where $\beta_z = 0$, to find that, with $b := (1-a)/2 > 0$,

$$\log\frac{\alpha_z}{c_\vartheta} = i\varphi(N_x - \mu_x) + O\Big(\varphi^2\mu_x^2 + \frac{1}{(\log x)^b}\Big), \tag{10.7}$$

since, by Lemma 4.4, $s_0(z)+1=i\varphi C+O\big(\varphi^2\big)$. Writing

$$E_z:=\frac{\log(\alpha_z/c_\vartheta)}{i\varphi},\quad R_x=N_x-\mu_x,$$

this implies

$$E_z=R_x+O(|\varphi|\mu_x^2+\frac{1}{|\varphi|(\log x)^b}\Big)\qquad(|\varphi|\leqslant\delta,\ x\geqslant 3). \tag{10.8}$$

We claim that (10.8) implies $E_z\ll 1$ and $R_x\ll 1$. The desired estimate (10.5) then follows from $E_z\ll 1$. To see this, apply (10.8) with $\varphi_j=1/(\log x)^{b/2}$, and then again with $\varphi_{j+1}=1/\big(\frac12 b\log_2 x\big)^3$, to get

$$E_{z_j}=E_{z_{j+1}}+O(1/|\log\varphi_j|),\quad \varphi_{j+1}=1/|\log\varphi_j|^3.$$

Starting with $\varphi_1=\varphi$, we iterate this equation to find

$$E_z=E_{z_k}+O\Big(\sum_{1\leqslant<k}\frac{1}{|\log\varphi_j|}\Big),$$

where $k$ is such that $\varphi_k$ lies between two suitable positive constants. Since $E_z\ll 1$ for $\varphi\asymp 1$, and $\sum_{1\leqslant j<k}1/|\log\varphi_j|\ll 1$, we have $E_z\ll 1$ for $|\varphi|\leqslant\delta$.

To see that $R_x\ll 1$, choose $\varphi=1/\mu_x^3$ in (10.8). □

## 11. Proof of Theorem 1.4

From the proof of Theorem 10.1 we know that $R_x\ll 1$, that is

$$N_x:=\frac{1}{B(x)}\sum_{n\in\mathcal{B}(x)}\nu(n)=\mu_x+O(1),$$

which coincides with (1.3). To prove (1.4), we will show that

$$M(x):=B(x)^{-1}\sum_{n\in\mathcal{B}(x)}\nu(n)^2=\mu_x^2+O(\mu_x), \tag{11.1}$$

with still $\mu_x=C\log_2 x$. This implies

$$\sum_{n\in\mathcal{B}(x)}(\nu(n)-\mu_x)^2\ll\mu_x B(x),$$

which plainly yields the desired result.

It remains to establish (11.1). Exponentiating (10.7), we have

$$\frac{\alpha_z}{c_\vartheta}=1+i\varphi(N_x-\mu_x)+O\Big(\varphi^2\mu_x^2+\frac{1}{(\log x)^b}\Big).$$

Since the real part of $\alpha_z/c_\vartheta-1$ is independent of $x$, we can write

$$\frac{\alpha_z}{c_\vartheta}=1+i\varphi(N_x-\mu_x)+J(\varphi)+iL(x,\varphi), \tag{11.2}$$

where $J(\varphi),L(x,\varrho)\in\mathbb{R}$ and $J(\varphi),L(x,\varphi)\ll\varphi^2\mu_x^2+1/(\log x)^b$. Considering $x$ such that $\varphi^2=1/(\log x)^b$ yields that $J(\varphi)\ll\varphi^2(\log\varphi)^2$.

Writing $z^{\nu(n)}=e^{i\varphi\nu(n)}=1+i\varphi\nu(n)-\frac12\varphi^2\nu(n)^2+O\big(\varphi^3\nu(n)^3\big)$, we obtain

$$B(x,z)=\frac{c_\vartheta x}{\log x}\left\{1+N_x i\varphi-M(x)\varphi^2/2+O\Big(|\varphi|^3\mu_x^3+\frac{1}{(\log x)^b}\Big)\right\}$$

for $\mu_x\varphi \ll 1$, by Lemmas 3.2 and 10.6. Now Theorem 10.1 and Lemma 4.4 yield

$$B(x,z) = \frac{\alpha_z x}{\log x}\left\{1 + \mu_x i\varphi - \tfrac{1}{2}\mu_x^2\varphi^2 + O\Big(\varphi^2\mu_x + \frac{1}{(\log x)^b}\Big)\right\}.$$

Equate the right-hand sides of the last two displays, divide by $c_\vartheta x/\log x$, insert the estimate (11.2) for $\alpha_z/c_\vartheta$, and compare the real parts of both sides: it only remain to select $\varphi = 1/\mu_x^2$ to obtain (11.1).

## 12. Proof of Theorem 1.3

Let $F_x(y)$ denote the left-hand side of (1.2) and let

$$\psi_x(\tau) := \frac{1}{B(x)}\sum_{n\in\mathcal{B}(x)} \mathrm{e}^{i\tau\{\nu(n)-C\log_2 x\}/\sqrt{V\log_2 x}} \qquad (\tau\in\mathbb{R}).$$

With $z = \mathrm{e}^{i\varphi}$, we have $z-1 = i\varphi - \frac{1}{2}\varphi^2 + O\big(\varphi^3\big)$ and $(z-1)^2 = -\varphi^2 + O\big(\varphi^3\big)$. Lemma 4.4 and Corollary 10.2 yield, uniformly for $x \geqslant 2$, $|\varphi| \leqslant \delta$,

$$B\big(x,\mathrm{e}^{i\varphi}\big)/B(x) = \{1+O(\varphi)\}(\log x)^{Ci\varphi - V\varphi^2/2 + O(\varphi^3)} + O\Big(\frac{1}{(\log x)^{1-a}}\Big).$$

Selecting

$$T := \delta\sqrt{V\log_2 x}, \quad \varphi := \frac{\tau}{\sqrt{V\log_2 x}},$$

we obtain

$$\psi_x(\tau) \ll (\log x)^{-V\varphi^2/3} = \mathrm{e}^{-\tau^2/3} \qquad (|\tau| \leqslant T),$$

provided $\delta$ is sufficiently small. From this point forward, the proof is identical to that of Theorem 1.2, with Lemma 10.6 replacing Lemma 6.4, and $a$ replacing $\varepsilon$.

## Acknowledgment

This work began when the second-named author visited the Institut Élie Cartan de Lorraine, whose support and hospitality are gratefully acknowledged.

Institut Élie Cartan, Université de Lorraine, B.P. 70239,
F–54506 Vandœuvre-lès-Nancy Cedex, France
*Email address*: `gerald.tenenbaum@univ-lorraine.fr`

Department of Mathematics, 351 West University Boulevard,
Southern Utah University, Cedar City, Utah 84720, USA
*Email address*: `weingartner@suu.edu`